%% file: formal.tex
\documentclass[10pt]{amsart}
\usepackage{amssymb,amsfonts,amsmath,amsthm}
\usepackage[all]{xy}
\input{header.tex}
\input{psfig.tex}

\begin{document}

\title{Formalism for Relative Gromov-Witten Invariants}
\author{Eric Katz}
\date{\today}
\address{Department of Mathematics, Duke University} \email{eekatz@math.duke.edu}
\begin{abstract}
  We develop a formalism for relative Gromov-Witten invariants
  following Li \cite{Li1,Li2} that is analogous to the Symplectic
  Field Theory of Eliashberg, Givental, and Hofer \cite{EGH}.
  This formalism allows us to express natural degeneration formulae in terms of generating
  functions and re-derive the formulae of Caporaso-Harris
  \cite{CH}, Ran \cite{Ra1}, and Vakil \cite{Va1}.  In addition, our framework gives a homology theory analogous to
SFT Homology.
\end{abstract}
\maketitle
    \section{Introduction}

\input{introformal.tex}
    \section{Background}
    \input{background.tex}
    \section{Generating Functions}

\input{generating.tex}
    \section{Degeneration Formulae}

\input{dgr.tex}
    \section{Examples}

\input{examples.tex}
    \section{Hamiltonian Formalism}

\input{hamiltonian.tex}
    \section{Localization Proof of Degeneration Formula}

\input{local.tex}
    \bibliographystyle{plain}
    \bibliography{formal}

\vspace{+10 pt} \noindent

\end{document}

%% file: header.tex
\newcommand{\proj}{\mathbb P}
\newcommand{\eff}{\mathbb F}
\newcommand{\aff}{\mathbb A}

\newcommand{\com}{\mathbb{C}}
\newcommand{\zee}{{\mathbb{Z}}}
\newcommand{\zeeg}{{\zee}_{\geq 0}}

\newcommand{\rat}{\mathbb{Q}}
\newcommand{\oh}{{\mathcal{O}}}
\newcommand{\ca}{{\mathcal{A}}}

\newcommand{\cf}{{\mathcal{F}}}

\newcommand{\cl}{{\mathcal{L}}}
\newcommand{\cz}{{\mathcal{Z}}}
\newcommand{\ct}{{\mathcal{T}}}
\newcommand{\cw}{{\mathcal{W}}}
\newcommand{\ex}{{\mathcal{X}}}
\newcommand{\cy}{{\mathcal{Y}}}

\newcommand{\cb}{{\mathcal{B}}}

\newcommand{\tcg}{{\widetilde{\mathcal{G}}}}
\newcommand{\ccr}{{\mathcal{R}}}
\newcommand{\ccp}{{\mathcal{P}}}

\newcommand{\ch}{{\mathcal{H}}}

\newcommand{\cm}{{\mathcal{M}}}
\newcommand{\cp}{{\mathcal{P}}}
\newcommand{\cq}{{\mathcal{Q}}}

\newcommand{\ma}{{\mathcal{MA}}}
\newcommand{\mz}{{\mathcal{MZ}}}
\newcommand{\my}{{\mathcal{MY}}}

\newcommand{\cmbar}{{\overline{\cm}}}
\newcommand{\acm}{{\mathfrak{M}}}
\newcommand{\mrn}{\cmbar_{0,n}}

\newcommand{\Top}{\operatorname{L^0}}
\newcommand{\Bot}{\operatorname{L^\infty}}
\newcommand{\Split}{\operatorname{Split}}
\newcommand{\Dil}{{\operatorname{Dil}}}
\def\Lnt#1{\operatorname{L}_{#1,\text{not top}}}
\def\Lnb#1{\operatorname{L}_{#1,\text{not bot}}}
\def\Le#1{\operatorname{L}_{#1,\text{ext}}}

\newcommand{\pt}{\text{pt}}
\newcommand{\ev}{\text{ev}}
\newcommand{\Ev}{\text{Ev}}
\newcommand{\Bl}{\text{Bl}}
\newcommand{\ft}{\text{ft}}
\newcommand{\id}{\text{id}}

\newcommand{\Sing}{\operatorname{Sing}}

\newcommand{\vdim}{{\operatorname{vdim\ }}}

\def\sig#1{\sigma_{#1}}
\def\vir#1{[#1]^{\text{vir}}}

\def\HS#1{{H_{\com^*}^{#1}}}
\newcommand{\tz}{{\tilde{z}}}
\newcommand{\thbar}{{\tilde{\hbar}}}
\newcommand{\tp}{{\tilde{p}}}
\newcommand{\tP}{{\tilde{P}}}
\newcommand{\tTheta}{{\tilde{\Gamma}}}
\newcommand{\Aut}{\operatorname{Aut}}

\newcommand{\coker}{\operatorname{coker}}
\newcommand{\im}{\operatorname{im}}

\def\nicechartonZ#1#2
{\xymatrix{
\ex \ar[d] \ar[r]^{#1} & Z[n]\ar[d]\\
S \ar[r]_{#2} & \aff^n}}
\def\nicechartonZmp#1#2
{\xymatrix{
\ex \ar[d] \ar[r]^{#1} & Z[n]\ar[d]\\
S \ar@<1ex>[u]^{\gamma}\ar[r]_{#2} & \aff^n}}
\def\nicechartonAmp#1#2
{\xymatrix{
\ex \ar[d] \ar[r]^{#1} & A[n]\ar[d]\\
S \ar@<1ex>[u]^{\gamma}\ar[r]_{#2} & \aff^n}}

\newcommand{\srarr}{\rightarrow}

\newcommand{\cs}{\com^*}

\theoremstyle{plain}
\newtheorem{theorem}{Theorem}
\newtheorem{lemma}[theorem]{Lemma}
\newtheorem{proposition}[theorem]{Proposition}
\newtheorem{corollary}[theorem]{Corollary}
\theoremstyle{definition}
\newtheorem{definition}[theorem]{Definition}

\newtheorem{example}[theorem]{Example}
\numberwithin{theorem}{subsection}

%% file: introformal.tex
Relative Gromov-Witten invariants following Li \cite{Li1,Li2} and
the Symplectic Field Theory of Eliashberg, Givental, and Hofer
\cite{EGH} are both theories of holomorphic curves with asymptotic
boundary conditions.  They have different sources: the theory of
relative Gromov-Witten invariants counts stable maps to a
projective manifold relative a divisor and is a systematization of
degeneration methods in enumerative geometry \cite{CH,Ra1,Va1};
Symplectic Field Theory, a generalization of Floer Homology. SFT
has an interesting formal structure involving a differential
graded algebra whose homology is an invariant of contact
structures.

In relative Gromov-Witten theory, one considers a pair $(Z,D)$
where $Z$ is a projective manifold and $D$ is a smooth, possibly
disconnected divisor in $Z$.  One looks at stable maps to $Z$
where all points of intersections of the map with $D$ are marked
and multiplicities at these points are specified.  To obtain a
proper moduli stack of such maps, one must allow the target to
degenerate to $\ _k Z=Z\sqcup_D P_1 \sqcup_D \dots \sqcup_D P_k$,
that is, $Z$ union a number of copies of $P=\proj_D (N_{D/Z}\oplus
1_D)$ the projective completion of the normal bundle to $D$ in
$Z$. Maps with a non-smooth target are said to be {\em split
maps}.  Li constructed a moduli stack of relative maps called
$\cm(\cz,\Gamma)$ for $\Gamma$, a certain kind of graph, and
constructed its virtual fundamental cycle.  This stack has an
evaluation map
$$\Ev_{\mz}:\cm(\cz,\Gamma)\srarr Z^m \times D^r$$
where $m$ and $r$ are the number of interior and boundary marked
points, respectively. {\em Relative Gromov-Witten invariants} are
given by evaluating pullbacks of cohomology classes by $\Ev$
against the virtual cycle.

It is natural to break the target $\ _k Z$ as the union of $\ _l
Z=Z\sqcup_D P_1 \sqcup_D \dots \sqcup_D P_l$ and $\ _{k-l-1}
P=P_{l+1}\sqcup_D \dots \sqcup_D P_k$. In fact, such splitting is
necessary to parameterize fixed loci in $\com^*$-localization in
the sense of \cite{K} and \cite{GP} in the relative framework
\cite{GV}. If we set $X=D$, and $L=N_{D/Z}$, the normal bundle to
$D$ in $Z$, one is led to study stable maps into the
projectivization of a line bundle $P=\proj_X(L\oplus 1_X)$
relative to the zero and infinity sections, $D_0$ and $D_\infty$
where two stable maps are declared equivalent if they can be
related by a $\com^*$-factor dilating the fibers of $P\srarr X$.
One can construct a moduli stack of such maps, $\cm(\ca,\Gamma)$
and its virtual cycle.  This moduli stack has certain natural line
bundles, called the target cotangent line bundles, $\Top$ and
$\Bot$ and has an evaluation map
$$\Ev_{\ma}:\cm(\ca,\Gamma)\srarr X^m \times X^{r_0} \times X^{r_\infty}$$
The {\em rubber invariants} are obtained by evaluating pullbacks
of cohomology by $\Ev$ map and powers of $c_1(\Bot)$ against the
virtual cycle.

The purpose of this paper is a systematic development of the
formal structure of the relative Gromov-Witten Invariants,
organized in generating functions.

We note here that the rubber invariants have been introduced
previously in the literature by Okounkov and Pandharipande
\cite{OPVir} and  by Graber and Vakil \cite{GV} as {\em maps to a
non-rigid target}.

In section 2, we recall the necessary background information to
describe the stacks $\cm(\cz,\Gamma_Z)$ and $\cm(\ca,\Gamma_A)$.
We show how to glue together such stacks to parameterize split
maps in a stack $\cm(\cz,\Gamma_Z * \Gamma_A)$. We describe
line-bundles on $\cm(\cz,\Gamma)$: $\Dil$ and $\Le{i}$; and
line-bundles on $\cm(\ca,\Gamma)$: $\Split$, $\Top$, $\Bot$,
$\Lnt{i}$, $\Lnb{i}$. These line-bundles have geometric meaning:
$\Le{i}$ is a line-bundle which has a section whose zero-stack
consists of maps $f:C\srarr\ _k Z$ so that the $i$th marked point
is not mapped to $Z\subset\ _k Z$ (counted with multiplicity);
$\Split$ is a line-bundle whose zero-stack is all split maps;
$\Lnt{i}$, where $i$ is the label of interior marked point, is a
line-bundle whose zero-stack consists of all split maps
$f:C\srarr\ _k P$ where $i$th marked point is not mapped to $P_k$;
$\Lnb{i}$ is its upside-down analog.

These line-bundles satisfy certain relations.  On
$\cm(\cz,\Gamma_Z)$:
\begin{eqnarray*}
\ev_i^*\oh(D)&=&\Le{i};
\end{eqnarray*}
and on $\cm(\ca,\Gamma_A)$:
\begin{eqnarray*}
\Top\otimes\Bot&=&\Split\\
\Top\otimes\ev_i^*L^\vee&=&\Lnt{i}\\
\Bot\otimes\ev_i^*L&=&\Lnb{i}.
\end{eqnarray*}

In section 3, we organize intersection numbers on
$\cm(\cz,\Gamma)$ and $\cm(\ca,\Gamma)$ into generating functions.
The intersection numbers on $\cm(\cz,\Gamma)$ of the form
$$\deg({\Ev_{\mz}}^* c \cap \vir{\cm(\cz,\Gamma)})$$
are organized into the {\em relative potential} $F$ which takes
values in a particular graded algebra $\cf$.  The intersection
numbers on $\cm(\ca,\Gamma)$,
$$\deg(c_1(\Bot)^l\cup {\Ev_{\ma}}^* c \cap \vir{\cm(\ca,\Gamma)})$$
are organized into the {\em rubber potential} $A$ which lies in an
algebra $\ccr$.  The algebra $\ccr$ acts on $\cf$ which
corresponds to joining curves in $\cm(\cz,\Gamma_Z)$ and
$\cm(\ca,\Gamma_A)$ to form split curves in
$\cm(\cz,\Gamma_Z*\Gamma_A)$.  Likewise, the multiplication
operation in $\ccr$ corresponds to joining curves in
$\cm(\ca,\Gamma_{Ab})$ to those in $\cm(\ca,\Gamma_{At})$.

In section 4, we prove {\em degeneration formulae} for the
relative and rubber potentials.  These degeneration formulae are
differential equations satisfied by the potentials and are
numerical consequences of the relations between line-bundles. Let
$F$ be the relative potential of a pair $(Z,D)$ and let
$A_{\lambda=0}$ be the rubber potential of the pair
$(D,L=N_{D/Z})$ without any powers of $c_1(\Bot)$.  Then, $F$
satisfies the differential equation

$$\sum_l N_{jl}\frac{\partial F}{\partial \theta_l}=\sum_l M_{jl}
\frac{\partial A_{\lambda=0}}{\partial \beta_l}\cdot F$$

\noindent where $\theta$ and $\beta$ are variables dual to
cohomology classes on $Z$ and $X$ respectively, $M_{jl}$ and
$N_{jl}$ are matrices that keep track of cohomology information,
and $\cdot$ is the action of $\ccr$ on $\cf$.

Given a pair $(X,L)$, the rubber potential satisfies the analogous
differential equation

$$\frac{\partial}{\partial \lambda}\frac{\partial A}{\partial
\beta_i}+\sum_j N_{ij}\frac{\partial A}{\partial
\beta_j}=\frac{\partial A_{\lambda=0}}{\partial \beta_i}*A$$

\noindent where $*$ is multiplication in $\ccr$.

In section 5, we work out several examples. We express the
rational rubber potential without powers of $c_1(\Bot)$ of
$(\proj^n,\oh(m))$ in terms of the Gromov-Witten invariants of
$\proj^n$ by a Kleiman-Bertini argument.  We use this rubber
potential to write down a degeneration formula for the relative
Gromov-Witten potential of $(\eff_n,D_\infty)$ and $(\proj^n,H)$
where $D_\infty\subset \eff_n$ is the infinity section of the
rational ruled surface of degree $n$, and $H$ is a hyperplane in
$\proj^n$.  This immediately yields the degeneration formulae of
of Caporaso-Harris \cite{CH}, Ran \cite{Ra1}, and Vakil
\cite{Va1}, phrased in the language of differential operators as
first stated by Getzler in \cite{G}.

In section 6, we construct a theory directly analogous to
Symplectic Field Theory.  One begin with a pair $(X,L)$ and
organizes a subset of the rubber invariants into a generating
function $H$ called the {\em Hamiltonian} that takes values in an
algebra $\ch$.

Given two interior marked points, one has the following formula
among divisors on $\cm(\ca,\Gamma)$:
\noindent
\begin{center}
\setlength{\unitlength}{.5cm}
\begin{picture}(17,5)
  \put(1,3)  {$\ev_2^*(c_1(L))-\ev_1^*(c_1(L))$}
  \put(8.5,3)  {$=$}
  \put(10.5,3.15){\line(1,-3){.5}}
  \put(10.5,2.85){\line(1, 3){.5}}
  \put(10.6,2.2) {\line(1,0) {.5}}
  \put(10.6,3.8) {\line(1,0) {.5}}
  \put(10.2,2)   {$2$}
  \put(10.2,3.6)   {$1$}
  \put(12.5,3)   {$-$}
  \put(14.5,3.15){\line(1,-3){.5}}
  \put(14.5,2.85){\line(1, 3){.5}}
  \put(14.6,2.2) {\line(1,0) {.5}}
  \put(14.6,3.8) {\line(1,0) {.5}}
  \put(14.2,2)   {$1$}
  \put(14.2,3.6)   {$2$}
\end{picture}
\end{center}
where the figures on the right specify certain loci of split
curves.   As a consequence of this formula, we have in $\ch$,
$$H^2=0$$

We can then define a differential on $\ch$ by the formula
$$D^H=Hf-(-1)^{\deg f}fH.$$
The homology of this complex, called Hamiltonain Homology is an
invariant of $(X,L)$ and is the algebraic geometric analog of the
SFT Homology of $S^1(L)$, the unit circle bundle of $L$.

In section 7, we give a direct proof of the degeneration formula
for the rubber potential using the technique of virtual
localization.

This paper draws most directly on the Relative Gromov-Witten
Invariants constructed by J. Li \cite{Li1,Li2} and the Symplectic
Field Theory of Eliashberg, Givental, and Hofer \cite{EGH}. Other
approaches to relative invariants include those of Gathmann
\cite{Gath1}, Ionel and Parker \cite{IP}, and A.-M. Li and Ruan
\cite{LiRuan}.

I would like to acknowledge the following for valuable
conversations: Y. Eliashberg, A. Gathmann, D. Hain, J. Li, and R.
Vakil. This paper, together with \cite{Ka2} is a revised version
of the author's Ph.D. thesis written under the direction of Y.
Elisahberg.

All varieties are over $\com$.

%\subsection{Projectivization Convention}

%%If $E$ is a vector bundle on a scheme $X$, let $\proj_X(E)$ denote
%%the projectization of $E$ where we use the old-fashioned geometric
%%notation where points in the projectivization represent lines.  In
%%this case, if $1_X$ is the trivial bundle on $X$, $\proj_X(E\oplus
%%1_X)$ is the projective completion of $E$.  The scheme
%%$$\proj_X(E\oplus 1_X)\setminus E$$
%%is called the infinity section.  In this sense, we follow
%%\cite{Fulton}.

%\end{document}

%% file: background.tex
%\input{preamble.tex}

%% refs for splitmaps extendedcompconent, topcomponent,
%% bottomcomponent

We discuss stacks of relative stable maps, $\mz=\cm(\cz,\Gamma)$
and stacks of maps to rubber, $\ma=\cm(\ca,\Gamma)$ where $\Gamma$
is a particular kind of graph. The material in this section is a
rephrasing of sections of \cite{Ka2}, some of which is
straightforward adaptation of \cite{Li1} and \cite{Li2}. While J.
Li does not construct $\ma$, our construction directly parallels
his.  We do change some notation from \cite{Li1} to suit our
purposes.

\subsection{Stacks of Relative Maps}

Consider a projective manifold $Z$ with a smooth divisor $D$.  Wew
review the construction of the stack of stable maps to $Z$
relative to $D$. Given an r-tuple of positive integers
$\mu=(\mu_1,\dots,\mu_r)$, consider a marked pre-stable curves
$$(C,x_1,\dots,x_m,p_1,\dots,p_r)$$
and maps
$$f:C\srarr Z$$
so that the divisor $f^*D$ is
$$f^*D=\sum_i \mu_i p_i.$$

To form a proper moduli stack of such maps, we must allow the
target to degenerate.  Let $L=N_{D/Z}$ be the normal bundle to $D$
in $Z$.  Let $P=\proj(L\oplus 1_D)$ be the projective completion
of $L$.  $P$ has two distinguished divisors, $D_0$ and $D_\infty$,
the zero and infinity sections of $L$.

\begin{definition} Let $_k Z$ be the union of $Z$ with $k$ copies
of $P$, $Z\sqcup_D P_1 \sqcup_D \dots \sqcup_D P_k$, the scheme
given by identifying $D\subset Z$ with $D_\infty \subset P_1$ and
$D_0\subset P_i$ with $D_\infty \subset P_{i+1}$ for
$i=0,1,\dots,k-1$.
\end{definition}

\begin{definition} Let $c:$ $_k Z\srarr Z$ be the {\em
collapsing map} that it the identity on $Z$ and projects each
$P_i$ to $D\subset Z$
\end{definition}

Note that $\Sing(_k Z)$, the singular locus of $_k Z$ is the
disjoint union of $k$ copies of $D$, which we label
$D_1,\dots,D_{k-1}$ where $D_i=D_\infty\subset P_i$.

\begin{definition} Let $\Aut(_k Z)=(\cs)^k$ be the group acting on
$_k Z$ where each factor of $\cs$ dilates the fibers of the
$\proj^1$ bundle $P_i\srarr X$.
\end{definition}

\begin{definition} Let $D\subset$ $_k Z$ denote the divisor
$D_0\subset P_k\subset$ $_k Z$.
\end{definition}

We need to specify the appropriate data for the moduli stack of
relative stable maps to $(Z,D)$.  Here we consider an algebraic
curve $C$ that is mapped to $_k Z$ by $f:C\srarr \ _k Z$ with
specified tangency to $D$. We must specify the topology of the
curve and the data of the marked points. There are two types of
marked points:
\begin{enumerate}
\item[(1)] {\em interior marked points} whose image under $f$ is
not mapped to $D$

\item[(2)] {\em boundary marked points} which are mapped to $D$ by
$f$.
\end{enumerate}
We will impose the condition that all points in $C$ mapped to $D$
will be marked.  The data of the curve is specified as follows.

\begin{definition} {A {\em relative graph} $\Gamma$ is the following
data:
\begin{enumerate}
\item[(1)] A finite set of vertices $V(\Gamma)$

\item[(2)] A genus assignment for each vertex
$$g:V(\Gamma)\srarr \zeeg$$

\item[(3)] A degree assignment for each vertex
$$d:V(\Gamma)\srarr B_1(Z)\equiv A_1(Z)/\tilde{ }_{\text{alg}}$$
that assigns the class of a curve modulo algebraic equivalence to
each vertex.

\item[(4)] A set $R=\{1,\dots,r\}$ labelling boundary marked
points together with a function assigning boundary marked points
to vertices
$$a_R:R\srarr V(\Gamma)$$

\item[(5)] A {\em multiplicity assignment} for each boundary
marked point
$$\mu:R\srarr \zee_{\geq 1}$$

\item[(6)] A set $M=\{1,\dots,m\}$ labelling interior marked
points together with an assignment to vertices
$$a_M:M\srarr V(\Gamma)$$
\end{enumerate}
}
\end{definition}

\begin{definition} Given two relative graphs $\Gamma$, $\Gamma'$
are said to be isomorphic if there is a bijection
$$q:V(\Gamma)\srarr V(\Gamma')$$
that commutes with the maps $g,d,a_R,a_M$.
\end{definition}

\begin{definition}
Let $\Gamma$ be a relative graph.  A {\em morphism to $_k Z$ of
type $\Gamma$} consists of a marked curve
$(C,x_1,\dots,x_{|M|},p_1,\dots,p_{|R|})$ and a morphism
$f:C\srarr\ _k Z$
\begin{enumerate}
\item[(1)] $C$ can be written as a disjoint union of pre-stable
curves $C_v$

\item[(2)] $C_v$ is a connected curve of arithmetic genus $g(v)$.

\item[(3)] The map
$$(c\circ f):C_v \srarr _k Z \srarr Z$$
has $(c\circ f)_*C_v=d(v)$.

\item[(4)] $x_i\in C_v$ for $v=a_M(i)$.  These are the interior
marked points. \label{interiormarkedpoint}

\item[(5)] $p_i\in C_v$ for $v=a_R(i)$.  These are the boundary
marked points.  \label{bdrymarkedpoint}

\item[(6)] $f^*D=\sum_{i\in R} \mu(i)p_i(s)$
\end{enumerate}
\end{definition}

\begin{definition} A morphism $f:C\srarr\ _k Z$ is said to be {\em
pre-deformable} if $f^{-1}(D_i)$ is the union of nodes so that for
$p\in f^{-1}(D_i)$ ($i=1,2,\dots,k)$, the two branches of the node
map to different irreducible component of $_k Z$ and that the
order of contact to $D_i$ are equal.
\end{definition}

\begin{definition} An {\em isomorphism} of morphisms $f,f'$ to $_k
Z$ consists of a diagram
%$$\begin{array}{ccccc}
%& & f & &\\
%& (C,x,p) & \srarr & _k Z &\\
%& h\downarrow & & \downarrow t &\\
%& (C',x',p') & \srarr & _k Z & \\
%& & f' & &
%\end{array}$$
$$\xymatrix{
(C,x,p) \ar[r]^>>>>>{f} \ar[d]_h & {\ _k Z} \ar[d]_t \\
(C',x',o') \ar[r]_>>>>>{f'} & {\ _k Z} }$$
 where $h$ is an isomorphism of
marked curves and $t\in \Aut(_k Z)$.
\end{definition}

\begin{definition} A morphism to $_k Z$ is said to be stable if it
has finitely many automorphisms.
\end{definition}

\begin{theorem} \cite{Li1} There is a Deligne-Mumford stack,
$\cm(\cz,\Gamma)$ parameterizing pre-deformable stable morphisms
to $_k Z$ for varying $k$.
\end{theorem}

In cases where it is understood, we will write $\mz$ for
$\cm(\cz,\Gamma)$

This moduli stack is constructed from a moduli functor by
considering stable maps to families of targets modelled on a
sequence of spaces and divisors $(Z[0],D[0])$,$(Z[1],D[1]),\dots$
defined inductively as follows
\begin{eqnarray*}
Z[0]&=&Z\\
D[0]&=&D\\
Z[n]&=&\Bl_{D[n-1]\times\{0\}} (Z[n-1]\times \aff^1)
\end{eqnarray*}
where $D[n]$ is the proper transform of $D[n-1]\times\aff^1$.
$Z[n]/Z$ possesses a $(\cs)^k$ group of automorphisms.  Note that
this construction is an iteration of deformation to the normal
cone.  $Z[n]$ posses a map to $\aff^n$.  Given a closed point
$x\in \aff^n$, the fiber over $x$ is $(Z[n])_x=\ _k Z$ where $k$
is the number of zeroes among $x$'s coordinates.

\begin{definition} \label{extendedcomponents} A map $f:C\srarr\ _k
Z$ is said to be {\em split} if $k\geq 1$.  The irreducible
components of $C$ that are mapped to $P_i\subset\ _k Z$ are said
to be {\em extended components}.
\end{definition}

\begin{definition} The {\em evaluation map} on
$\mz=\cm(\cz,\Gamma)$ is a map
$$\Ev:\mz\srarr Z^m\times D^r$$ given on a relative stable map
$(C,f)$ by
$$(x_1,\dots,x_m,p_1,\dots,p_r)\hookrightarrow C\srarr \ _k Z \srarr
Z.$$
\end{definition}

We will write $\ev_i:\mz\srarr Z$ or $\ev_i:\mz\srarr D$ to denote
the evaluation map at one of the interior or boundary marked
point.

\begin{theorem} \cite{Li2} $\mz$ carries a virtual cycle of
complex dimension
\begin{eqnarray*}
\vdim \mz&=&\sum_{v\in V(\Gamma)}(\dim
Z-3)(1-g(v))\\
&&+<c_1(TZ)-D,d(v)>)+|R|+|M|.
\end{eqnarray*}
\end{theorem}

\subsection{Stack of Maps to Rubber}

In constructing $\cm(\cz,\Gamma)$, we had to consider stable maps
to $_k Z$ which was $Z$ union a chain of $P$'s.  It is useful to
consider stable maps to the chain of $P$'s subject to
automorphisms.  We call these maps to rubber.

Let $X$ be a projective manifold and $L$ a line bundle on $X$. Let
$P=\proj_X(L\oplus 1_X)$, and let $X_0$ and $X_\infty$ denote the
zero and infinity sections.  We study stable maps to $P$ relative
to $X_0$ and $X_\infty$ where we mod out by a $\cs$-factor that
dilates the fibers.  Again, the target $P$ may degenerate.

\begin{definition} Let $_k P$ be the union of $k+1$ copies of $P$,
$$_k P =P_0 \sqcup_X P_1 \sqcup_X \dots \sqcup_X P_k$$
gluing $X_0\subset P_i$ to $X_\infty\subset P_{i+1}$ for
$i=0,\dots,k-1$.
\end{definition}

$_k P$ has distinguished divisors $D_\infty=X_\infty\subset P_0$
and $D_0=X_0 \subset P_k$.

\begin{definition} Let $\Aut(_k P)=(\cs)^{k+1}$ act on $_k P$ by
dilating fibers of $P\srarr X$.
\end{definition}

\begin{definition} A {\em rubber graph} $\Gamma$ is the following
data:
\begin{enumerate}
\item[(1)] a finite collection of vertices $V(\Gamma)$

\item[(2)] A genus assignment for each vertex
$$g:V(\Gamma)\srarr \zeeg.$$

\item[(3)] A degree assignment for each vertex
$$d:V(\Gamma)\srarr B_1(X)=A_1(X)/\tilde{ }_{\text{alg}}.$$

\item[(4)] Sets $R_0=\{1,\dots,r_0\}$,
$R_\infty=\{1,\dots,r_\infty\}$ labelling boundary marked points
together with a function assigning boundary marked points to
vertices
\begin{eqnarray*}
a_0&:&R_0\srarr V(\Gamma)\\
a_\infty&:&R_\infty\srarr V(\Gamma)
\end{eqnarray*}

\item[(5)] A {\em multiplicity assignment} for boundary marked
points
\begin{eqnarray*}
\mu^0&:&R_0\srarr \zee_{\geq 1}\\
\mu^\infty&:&R_\infty\srarr \zee_{\geq 1}
\end{eqnarray*}

\item[(6)] A set $M=\{1,\dots,m\}$ labelling interior marked
points together with an assignment to vertices
$$a_M:M\srarr V(\Gamma)$$
\end{enumerate}
\end{definition}

Definitions of morphisms to $_k P$ are analogous to morphisms to
$_k Z$ with $Z$'s replaced with $P$'s and the following
modifications. The degree assignment is
$$d(v)\in B_1(X)$$
We have marked points
$$(x_1,\dots,x_m,p^0_1,\dots,p^0_{|R_0|},p^\infty_1,\dots,p^\infty_{|R_\infty|})\subset
C$$ so that

\begin{eqnarray*}
f^*D_0&=&\sum_{i\in R_0} \mu^0(i)p^0_i\\
f^*D_\infty&=&\sum_{i\in R_\infty} \mu^\infty(i) p^\infty_i
\end{eqnarray*}

The {\em multiplicity condition} relates the multiplicities to
$D_0$ and $D_\infty$, to the degree:

\begin{lemma} \label{multcond} If $\cm(\ca,\Gamma)$ is nonempty
then for each vertex $v\in\Gamma$ we have
$$\sum_{p\in a_0^{-1}(v)} \mu^0(v)-\sum_{p\in a_\infty^{-1}(v)}
\mu^\infty(v)=<c_1(L),d(v)>$$
\end{lemma}

\begin{proof} One uses $X_0=X_\infty+\pi^*c_1(L)$ for each copy of $P$ in
the target.
\end{proof}

\begin{theorem} For a rubber graph $\Gamma$, $\ma=\cm(\ca,\Gamma)$
is a proper Deligne-Mumford stack.
\end{theorem}

\begin{theorem} $\ma$ carries a virtual cycle of complex dimension
\begin{eqnarray*}
\vdim \ma&=&\sum_{v\in V(\Gamma)}((\dim
X-2)(1-g(v))\\
&&+<c_1(TX),d(v)>)+|R_0|+|R_\infty|+|M|-1.
\end{eqnarray*}
\end{theorem}

We have analogous evaluation maps $\ev_i$ at the interior and
boundary marked points (mapping to $D_0$ and $D_\infty$).

\begin{definition} {The {\em evaluation map} on
$\ma=\cm(\ca,\Gamma)$ is
$$\Ev:\ma\srarr X^n\times X^{|R_0|}\times X^{|R_\infty|}$$}
\end{definition}

This moduli stack is constructed from a moduli functor by
considering stable maps to families of targets modelled on a
sequence of spaces and divisors
$(A[0],D_0[0],D_\infty[0])$,$(A[1],D_0[1],D_\infty[1]),\dots$
where
\begin{eqnarray*}
A[0]&=&P\\
D_0[0]&=&X_0\\
D_\infty[0]&=&X_\infty\\
A[n]&=&\Bl_{D_0[n-1]\times\{0\}} (A[n-1]\times \aff^1)
\end{eqnarray*}
where $D_0[n]$ is the proper transform of $D_0[n-1]\times\aff^1$.
and $D_\infty[n]$ is the inverse image of
$D_\infty[n-1]\times\aff^1$. $A[n]/X$ possesses a $(\cs)^{k+1}$
group of automorphisms

\begin{definition} \label{splitmaps} A {\em split map} in $\ma$ is
a map $f:C\srarr\ _k P$ where $k\geq 1$, that is, the target is
not smooth.
\end{definition}

\begin{definition} \label{topcomponent}
For a map $f:C\srarr\ _k P$ in $\ma$, the irreducible components
of $C$ that are mapped to $P_k$ are said to be the {\em top
components} while the components of $C$ that are mapped to $P_0$
are said to be the {\em bottom components}.
\end{definition}

We should explain our top/bottom convention. In $Z$, moving
towards $D$ is considered moving towards the top. In $P$, $D_0$ is
considered the top while $D_\infty$ is the bottom. This slightly
odd convention makes sense in that the most natural choice for
$(X,L)$ is $(D,N_{D/Z})$.  In this case, the zero section of $P$
is identified with $D$ and the normal bundle to $D_0$ in $P$ is
equal to the normal bundle to $D$ in $Z$. Therefore, $D_0\subset
P$ like $D\subset Z$ is on top.

\subsection{Trivial Cylinders}

We will single out certain connected components of curves
parameterized by $\ma$.  These are the so called trivial cylinders
which will be significant when we encode the data of the moduli
space into generating functions.

\begin{definition} \label{trivialcyl} Let $\Gamma$ be a rubber
graph. A vertex $v$ is said to correspond to a {\em trivial
cylinder} of degree $r$ if
\begin{enumerate}
\item[(1)] $g(v)=0$.

\item[(2)] $d(v)=0$

\item[(3)] $a_0^{-1}(v)$ is a single point.

\item[(4)] $a_\infty^{-1}(v)$ is a single point.

\item[(5)] $\mu^0(a_0^{-1}(v))=\mu^\infty(a_\infty^{-1}(v))=r$.

\item[(6)] $A_M^{-1}(v)$ is empty.
\end{enumerate}
\end{definition}

A trivial cylinder corresponds to a connected component of a map
to rubber.  This map is from a chain of $k$ $\proj^1$'s to $\ _k
P$ where each $\proj^1$ is mapped to a fiber of $P_i\srarr X$ and
is of degree $r$ and totally ramified at $X_0$ and $X_\infty$.

Note that if $\Gamma$ has a single vertex corresponding to a
trivial cylinder, then there are no rubber maps of type $\Gamma$
that are not invariant under the $\com^*$-action that dilates the
fibers of $P$.  Therefore, there are no stable rubber maps and the
moduli space is empty. This does not rule out morphisms of type
$\Gamma$ which has a component which is a trivial cylinder.  In
fact, one can add a trivial cylinder component to any family.

\subsection{Gluing Moduli Stacks}

Consider a projective manifold $Z$, together with a smooth divisor
$D$.  We will consider a relative moduli stack $\cm(\cz,\Gamma_Z)$
corresponding to $(Z,D)$ and a rubber moduli stack
$\cm(\ca,\Gamma_A)$ corresponding to $(X=D,L=N_{D/Z})$ where
$N_{D/Z}$ is the normal bundle to $D$ in $Z$.  One can join maps
in $\mz$ to maps in $\ma$ if certain conditions are met. Likewise
under particular conditions can join maps in
$\ma_1=\cm(\ca,\Gamma_1)$ to maps in $\ma_2=\cm(\ca,\Gamma_2)$. We
make these conditions precise below.

\begin{definition} Let
$\Gamma_Z$ be a relative graph and $\Gamma_A$ be a rubber graph.
Suppose that $L:RZ\srarr RA_\infty$ is a bijection from the
labelling sets for boundary marked points in $\Gamma_Z$ to the
labelling sets for boundary marked points mapping to $D_\infty$ in
$\Gamma_A$ so that
$$\mu_Z(q)=\mu^\infty_A(L(q)).$$
Let
$$J:M_Z\sqcup M_A\srarr \{1,\dots,|M_Z|+|M_A|\}$$
be a bijection between the labelling sets of the interior marked
points and a set of $|M_Z|+|M_A|$ elements.  We call the data
$(\Gamma_A,\Gamma_Z,L,J)$ a {\em graph join quadruple}.
\end{definition}
Colloquially, we've matched boundary marked points on $\Gamma_Z$
and $\Gamma_A$ with the same multiplicity.

\begin{definition} \label{graphjoin}
Define the {\em graph join} $\Gamma_A*_{L,J}\Gamma_Z$ to be the
following relative graph.  Let the graph $\Delta$ be obtained by
taking as vertices the vertices of $\Gamma_Z$ and $\Gamma_A$ and
for every $q\in RZ$, place an edge between the vertices
corresponding to $q$ and $L(q)$. Let $\Gamma_A*_{L,J}\Gamma_Z$ be
given as follows.  The vertices of
$\Gamma=\Gamma_A*_{L,J}\Gamma_Z$ are the connected components of
$\Delta$. Let $b_Z:V(\Gamma_Z)\srarr V(\Gamma)$, and
$b_A:V(\Gamma_A)\srarr V(\Gamma)$ be the functions taking vertices
of $\Gamma_A$ and $\Gamma_A$ to the components in $\Delta$
containing them. For $v\in V(\Gamma)$, let $\Delta_v$ be the
connected component of $\Delta$ corresponding to $v$.  Define the
data for $\Gamma_A*_{L,J}\Gamma_Z$ as follows:
\begin{enumerate}
\item[(1)] $g(v)=(\sum_{w\in b_Z^{-1}(v)} g(w))+(\sum_{w\in
b_A^{-1}(v)} g(w))+\dim(h^1(\Delta_v))$

\item[(2)] $d(v)=(\sum_{w\in b_Z^{-1}(v)} d(w))+(\sum_{w\in
b_A^{-1}(v)} i_*d(w))$ where $i:X\srarr Z$ is the inclusion and
$i_*:B_1(X)\srarr B_1(Z)$ is the induced map.

\item[(3)] $R=RA_0$ with $a_R:R\srarr V(\Gamma)$ given by
$$a_R=b_A \circ a_0$$

\item[(4)] $\mu:R\srarr \zee_{\geq 1}$ given by
$$\mu_R=\mu^0$$

\item[(5)] $M=\{1,\dots,|M_Z|+|M_A|\}$ with assignment function
$a:M\srarr V(\Gamma)$ given for $k\in J(M_Z)$ by
$$a(k)=b_Z \circ a_{MZ}\circ J^{-1}$$
while for $k\in J(M_A)$ by
$$a(k)=b_A \circ a_{MA}\circ J^{-1}$$
\end{enumerate}
\end{definition}

Given $(\Gamma_Z,\Gamma_A,L,J)$ as above, consider the evaluation
map at the boundary marked points on $\cm(\cz,\Gamma_Z)$ followed
by a map $L_*:D^r\srarr D^r$ which reorders the products of $D^r$
according to $L$:
$$L_*\circ Ev_R:\cm(\cz,\Gamma_Z)\srarr D^r\srarr D^r$$
and the evaluation map at the boundary marked points mapping to
$D_\infty\cong X$ on $\cm(\ca,\Gamma_A)$,
$$\Ev_{R_\infty}:\cm(\ca,\Gamma_A)\srarr D^r.$$

\begin{theorem} \cite{Li1} There is a morphism
$$\Phi_{\Gamma_Z,\Gamma_A,L,J}:\cm(\ca,\Gamma_A)\times_{D^r} \cm(\cz,\Gamma_Z)\srarr \cm(\cz,\Gamma_A*_{L,J}\Gamma_Z).$$
\end{theorem}

\begin{definition} \label{zuniona} {Let the stack
$\cm(\ca\sqcup\cz,\Gamma_A\sqcup_{L,J} \Gamma_Z)$ be the image
stack of $\Phi$ in $\cm(\cz,\Gamma_A*_{L,J} \Gamma_Z)$.}
\end{definition}

\begin{definition} An
{\em automorphism of $RZ$} is a permutation
$$\sigma:RZ\srarr RZ$$
so that $\mu_Z(\sigma(i))=\mu_Z(i)$ and
$a_{RZ}(\sigma(i))=a_{RZ}(i)$.  The group of all such
automorphisms is denoted by $\Aut_{\Gamma_Z}(RZ)$.  Likewise, we
define $\Aut_{\Gamma_A}(RA_0)$ and $\Aut_{\Gamma_A}(RA_\infty)$.

Given $L:RZ\srarr RA_\infty$, we may define
$\Aut_{\Gamma_A,\Gamma_Z,L}(RZ,RA_\infty)$ as the subgroup of
$\Aut_{\Gamma_Z}(RZ)\times\Aut_{\Gamma_A}(RA_\infty)$ such that
for $(\sigma,\tau)\in
\Aut_{\Gamma_Z}(RZ)\times\Aut_{\Gamma_A}(RA_\infty)$ we have
$L(\sigma(i))=\tau(L(i))$ for $1\leq i\leq |RZ|$.
\end{definition}

\begin{lemma} (\cite{Li1}, Prop 4.13) $\Phi$ is finite and \'{e}tale onto its image of degree equal
to
$$|\Aut_{\Gamma_A,\Gamma_Z,L}(RZ,RA_\infty)|$$ at every integral
substack of $\cm(\ca\sqcup\cz,\Gamma_A \sqcup_{L,J} \Gamma_Z)$.
\end{lemma}

\begin{definition} \label{multgraphjoin} Let
$\Upsilon=(\Gamma_Z,\Gamma_A,L,J)$ be a graph join quadruple Let
the {\em boundary multiplicity} $m(\Upsilon)$ be given by
$$m(\Upsilon)=\prod_{i\in RZ} \mu_Z(i).$$
\end{definition}

\begin{definition} \label{joinclasses} {Two quadruples are said to
be {\em join-equivalent} if they give the same image under $\Phi$}
\end{definition}

\begin{proposition} \label{totalgluing} Consider a join-equivalence
class of quadruples
$$[\Upsilon]=[(\Gamma_A,\Gamma_Z,L,J)].$$
Let $N=\cm(\ca\sqcup\cz,\Upsilon)$. Let
$$M_\Upsilon=\coprod_{(\Gamma_Z',\Gamma_A',L',J')}
\cm(\ca,\Gamma_A')\times_{D^r}\cm(\cz,\Gamma_Z')$$ where the
disjoint union is over quadruples join-equivalent to $\Upsilon$.
Then $\Phi_{[\Upsilon]}:M\srarr N$ is an \'{e}tale map of degree
$$|M_Z|!|M_A|!(|RZ|!)^2$$
\end{proposition}

\begin{proof}
This follows from the previous lemma and that there are
$$|M_Z|!|M_A|!\frac{(|RZ|!)^2}{|\Aut_{\Gamma_A,\Gamma_Z,L}(RZ,RA_\infty)|}$$
elements in $(\Gamma_Z,\Gamma_A,L,J)$'s graph join-equivalence
class.
\end{proof}

Likewise, we may define graph-join for rubber graphs, $\Gamma_t$,
$\Gamma_b$ (where $t$ and $b$ stand for top and bottom).  Let
$L:R_{b0}\srarr R_{t\infty}$ be a bijective function satisfying
$$\mu_b^0(q)=\mu_t^\infty(L(q)).$$
Let
$$J:M_b\sqcup M_t\srarr \{1,\dots,|M_b|+|M_t|\}$$
be a bijective map.  Then we define the {\em graph join}, a rubber
graph $\Gamma=\Gamma_t*_{L,J} \Gamma_b$ as above, except that
instead of condition (3) above, we have
$$R_0=R_{t0},\ \ a_0=b_{A_t}\circ a_{0t}$$
$$\mu^0=\mu_t^0$$
$$R_\infty=R_{b\infty},\ \ a_\infty=b_{A_b}\circ a_{\infty}b$$
$$\mu^\infty=\mu_b^\infty.$$

\begin{definition} Let $\Upsilon=(\Gamma_{A_t},\Gamma_{A_b},L,J)$ be a
quadruple describing a decomposition in
$\cm(\ca,\Gamma_{A_t}*_{L,J}\Gamma_{A_b})$. Define $m(\Upsilon)$
by
$$m(\Upsilon)=\prod_{i\in RA_t} \mu_t^\infty(i).$$
\end{definition}

Now, let $r=|R_{b0}|=|R_{t\infty}|$.  Exactly as above, we have

\begin{theorem} \cite{Li1}{There is a morphism
$$\Phi:\cm(\ca,\Gamma_{A_t})\times_{D^r} \cm(\ca,\Gamma_{A_b})\srarr \cm(\ca,\Gamma_{A_t}*_{L,J}\Gamma_{A_b}),$$
\'{e}tale of degree
$|\Aut_{\Gamma_{A_b},\Gamma_{A_t},L}(RA_{b0},RA_{t\infty})|$ }
\end{theorem}

\begin{corollary} {Consider a moduli stack $N=\cm(\ca,\Gamma_{A_b}
\sqcup_{L,J} \Gamma_{A_t})$.  Let
$$M=\coprod_{(\Gamma_{A_t}',\Gamma_A',L',J')}
\cm(\ca,\Gamma_{A_b}')\times_{D^r}\cm(\ca,\Gamma_{A_t}')$$ where
the disjoint union over $(\Gamma_{At},\Gamma{Ab},L,J)$'s
join-equivalence class. The $M\srarr N$ is an \'{e}tale map of
degree
$$|M_{A_b}|!|M_{A_t}|!(|R_{A_b0}|!)^2$$}
\end{corollary}

\subsection{Line Bundles on Moduli Stacks}

The moduli stacks carry line-bundles with particular geometric
meaning.

Given a relative graph $\Gamma_Z$, $\cm(\cz,\Gamma_Z)$ has
canonically defined line-bundles
\begin{enumerate}
\item[(1)] $\Dil$, a line bundle that has a section whose zero
stack is supported on all split curves.

\item[(2)] $\Le{i}$ where $i$ is a distinguished interior marked
point, a line bundle that has a section whose zero stack is
supported on split curves where $i$ lies on an extended component.
\end{enumerate}

It will be shown that $c_1(\Dil)$ on $\mz$ is (counted with
multiplicity) the locus of split maps and $c_1(\Le{i})$ is a
weighted count of split maps with $i$ on an extended component.

For $\Gamma_A$, $\cm(\ca,\Gamma_A)$ has the following line-bundles
\begin{enumerate}
\item[(1)] $\Top$, the target cotangent line-bundle at $X_0$.
$c_1(\Top)=\Psi_0$, the target $\Psi$ class of \cite{FP}.

\item[(2)] $\Bot$, the target cotangent line-bundle at $X_\infty$
which is $\Top$'s upside-down analog. $c_1(\Bot)=\Psi_\infty$.

\item[(3)] $\Split$, the Split bundle which has a section whose
zero stack is supported on all split maps (Definition
\ref{splitmaps}).

\item[(4)] $\Lnt{i}$, the not-top bundle with respect to a
distinguished interior marked point $i$.  This bundle has a
section whose zero stack is supported on split maps where the
$i$th marked point is not on a top component.

\item[(5)] $\Lnb{i}$, the not-bottom bundle with respect to a
interior marked point $i$.  This bundle has a section whose zero
stack is supported on split maps where the $i$th marked point is
not on a bottom component.
\end{enumerate}

$\Top$, which is defined in terms of an atlas, has the following
intuitive description: given a map to rubber, $(C,f)$ in $\ma$,
consider the target of $f$, $\ _k P$ which has a top component
$P_k$.   Let $\widehat{C}$ be the component of $C$ mapping to
$P_k$. There is a $\cs$ family of maps
$\widehat{f}:\widehat{C}\srarr P$ that occur as the restriction of
f; these $\cs$ families fit together to give a $\cs$ bundle; the
associated $\com$ bundle is $\Top$.  $\Bot$ is the analogous
bundle where we consider the bottom component.

$\Top$ and $\Bot$ can be given an interpretation in the stack of
rational sausages, the substack of $\acm_{0,2}$ consisting of
pre-stable curves so that the two marked points lie on different
sides of every node. $\Top$ and $\Bot$ are equal to the pullbacks
of the cotangent line classes at the two marked points. See
\cite{GV} for an elaboration.

\begin{theorem} \cite{Ka2}\label{mpiobb} The line-bundles satisfy the following relations:
On $\cm(\cz,\Gamma_Z)$,
\begin{eqnarray*}
\ev_i^*\oh(D)&=&\Le{i}\\
\end{eqnarray*}
and on $\cm(\ca,\Gamma_A)$,
\begin{eqnarray}
\Top\otimes\Bot&=&\Split\\
\Top\otimes\ev_i^*L^\vee&=&\Lnt{i}\\
\Bot\otimes\ev_i^*L&=&\Lnb{i}
\end{eqnarray}
\end{theorem}

If we consider the stack of rational sausages where $\Top$ and
$\Bot$ are the restriction of $\psi$ classes on $\acm_{0,2}$, then
the (1) is the pullback of the genus 0 recursion relation of Lee
and Pandharipande \cite{LP}.  A proof of the enumerative
consequences of (2) and (3) is given in section 7.

%\bibliographystyle{plain}
%\bibliography{thesis}

%\end{document}

%% file: generating.tex
An important ideas in Gromov-Witten theory, originating
\cite{Wit}, is that of organizing invariants in generating
functions. Relations satisfied by the invariants become
differential equations for the generating function. In this
section, we define generating functions for relative and rubber
invariants motivated by Symplectic Field Theory \cite{EGH}.

%%The generating functions take values in a particular algebra. This
%%algebra is graded with the grading coming from that of the
%%cohomology of the target of the evalatuation map.  In addition,
%%there is a Weyl algebra structure. This Weyl algebra structure
%%will give a particularly nice way of expressing the relations
%%among invariants.  The definitions in this section are motivated
%%by Symplectic Field Theory \cite{EGH}.

\subsection{Relative Potential}

Let us consider a pair $(Z,D)$ where $Z$ is a projective manifold
and $D\subset Z$ is a smooth divisor on $Z$.  The generating
function of the relative invariants takes values in a particular
graded algebra.

%%\begin{definition} Let $B_1(Z)=A_1(Z)/{\tilde{\ }}_\text{alg}$ be
%%the algebraic one-cycles up to algebraic equivalence.
%%\end{definition}

Let us specify the following data: an Euler characteristic,
$\chi$; a curve class (up to algebraic equivalence), $d\in
B_1(Z)$; a number of interior marked points: $m$; a number of
boundary marked points, $r$; and a r-tuple of multiplicities to
$D$, $(s_1,s_2,\dots,s_r)$. Consider the set $\Xi$ of relative
graphs $\Gamma$ so that
\begin{enumerate}
\item[(1)] $\sum_v d(v)=d$

\item[(2)] $\sum_v (2-2g(v))=\chi$

\item[(3)] $|M|=m$

\item[(4)] $|R|=r$

\item[(5)] $(\mu(1),\mu(2),\dots,\mu(r))=(s_1,s_2,\dots,s_r)$
\end{enumerate}
To each $\Gamma$, we associate the moduli space
$\mz=\cm(\cz,\Gamma)$ which has an evaluation map at the interior
and boundary marked points,
$$\Ev:\mz\srarr (Z)^m\times (D)^r.$$

\begin{definition} Given cohomology classes
$$e_1,\dots,e_n\in
H^*(Z)$$ and $$c_1,\dots,c_r\in H^*(D),$$ define the {\em
correlator} by
$$<e_1,\dots,e_n\cdot c_1,\dots,c_r>_{\chi,A,(s_1,\dots,s_r)}$$
$$=\sum_{\Gamma\in\Xi} \Ev^*(e_1\times\dots\times e_n\times
c_1\times\dots\times c_r)\cap \vir{\cm(\cz,\Gamma)}.$$
\end{definition}

By stability considerations, the above sum over relative graphs is
finite.
%%\begin{lemma} The above sum over relative graphs is
%%finite.\end{lemma}

%%\begin{proof} There are finitely many relative graphs in $\Xi$ that
%%correspond to a non-empty moduli space.  First note that there are
%%finitely many vertices with $g(v)=d(v)=0$ since each such vertex
%%must have at least three marked points assigned to it. Likewise,
%%there are only finitely many vertices with $d(v)=0$ and $g(v)=1$.
%%Since every vertex with $d(v)=0$, $g(v)\geq 2$ gives a negative
%%contribution to the sum $\sum_v (2-2g(v))=\chi$, there are only
%%finitely many vertices with $d(v)=0$.  Since $d(v)$ must be in the
%%effective cone, and $\sum_v d(v)=d$, there are only finitely many
%%possibilities for degrees of the non-contracted components.
%%\end{proof}

Let $e_1,e_2,\dots,e_l\in H^*(Z)$ be a homogeneous basis of
$H^*(Z)$.  Let $c_1,c_2,\dots,c_k$ be a homogeneous basis of
$H^*(D)$.  Let $\rat[B_1(Z)]$ be the group algebra on $B_1(Z)$,
generated as a vector space by elements of $B_1(Z)$, equipped with
multiplication
$$\tz^d_1\cdot \tz^d_2=\tz^{d_1+d_2}$$
where $d_1,d_2\in B_1(Z)$.

Consider the graded super-commutative algebra over $\rat[B_1(Z)]$
freely generated by $\thbar^{-1}$, $\thbar$,
$\theta_1,\theta_2,\dots,\theta_l$, and, for every positive
integer $n$, elements $\tp_{n,1},\tp_{n,2},\dots,\tp_{n,k}$ with
the following degrees
\begin{eqnarray*}
\deg \thbar&=&-2(\dim Z-3)\\
\deg \tz^d&=&2<c_1(TX),d>\\
\deg \theta_i&=&2-\deg e_i\\
\deg \tp_{n,j}&=&2-\deg c_i-2j\\
\end{eqnarray*}
$\hbar$ will be a formal variable corresponding to one half of the
Euler characteristic, $\tz^d$ to degree, $\theta_i$ to interior
marked points that are mapped to a cycle Poincare-dual to $e_i$,
and $\tp_{n,j}$ to boundary marked points with multiplicity $n$
and mapped to a cycle on $D$ Poincare-dual to $c_j$.

We define $\cf$ to be a partial completion of the above algebra.
We look at Laurent series in $\thbar$ whose coefficients are
polynomials in the $\tp$-variables whose coefficients are power
series in the $\theta$ variables.

Let $\tcg$ be the noncommutative algebra of power series in
$\theta_i e_i$. Let $\cp$ be the noncommutative algebra of power
series in $\tp_{n,j}c_j$.

Define the $\cf$-correlator to be
$$<\theta_{i_1}e_{i_1}\ \dots\ \theta_{i_n}e_{i_n},\tp_{n_1,j_1}c_{j_1}\ \dots\ \tp_{n_r,j_r}c_{j_r}>_{\chi,A}$$
$$=\theta_{i_1}\dots\theta_{i_n}\tp_{n_1,j_1}\dots\tp_{n_r,j_r}<e_{i_1},\dots,e_{i_n}\cdot
c_{j_1}\dots c_{j_r}>_{\chi,A,(n_1,n_2,\dots,n_r)}$$

Extend the $\cf$-correlator multi-linearly to a map
$$(,)_{g,A}:\tcg\otimes\cp\srarr\cf$$

Let $\tTheta\in\tcg$, $\tP\in\cp$ be given by
\begin{eqnarray*}
\tTheta&=&\sum_{l\geq 0} \frac{1}{l!} \left(\sum \theta_i
e_i\right)^l,\\
\tP&=&\sum_n \frac{1}{n!}\left(\sum_{k,i} \tp_{k,i}c_i\right)^n.
\end{eqnarray*}

\begin{definition} The {\em relative potential} of $(Z,D)$ is
$F\in \cf$ defined by
$$F=\sum_{g\geq 0} \sum_{d\in B_1(Z)} <\tTheta,\tP>_{g,d}
\hbar^{g-1} z^d.$$
\end{definition}

Note that $F$ is indeed in $\cf$.  A given coefficient of $z^A$ is
a polynomial in the $\tp$ variables because the number of $\tp$
variables is less than $A\cdot D$.  By our choice of degrees for
the formal variable and by the virtual dimension of $\mz$, $F$ is
a homogeneous element of degree of $0$.  Our usage of $F$
disagrees with that of \cite{EGH} because we consider disconnected
stable maps.

\begin{example} For the target $(\proj^1,\infty)$, the relative
potential is
$$F=\exp(\hbar^{-1}(\frac{1}{2!}\theta_0^2\theta_1+e^{\theta_1}p_{1,1}z)+h^0(-\frac{1}{24}\theta_1))$$
where $\theta_0$ and $\theta_1$ are dual to the classes
$[\proj^1], [\pt]$.  The $\hbar^{-1}\frac{1}{2!}\theta_0^2$
corresponds to contracted rational curves,
$\hbar^{-1}e^{\theta_1}p_{1,1}z$ to degree one rational curves
with arbitrarily many marked points, and
$-\hbar^0\frac{1}{24}\theta_1$ to contracted elliptic curves.
\end{example}

\begin{example} For the target $(\proj^2,L)$, the relative
potential is
$$F=\exp(F_{d=0}+F_{d\geq 0})$$
where $F_{d=0}$ and $F_{d\geq 0}$ correspond to connected degree
$0$ and to positive degree maps, respectively:
\begin{eqnarray*}
F_{d=0}&=&-\hbar^0\frac{1}{8}\theta_1+\hbar^0(\frac{1}{2!}\theta_2\theta_0^2+\frac{1}{2!}\theta_1^2\theta_0)\\
F_{d\geq 0}&=&
\hbar^{-1}\theta_2p_{1,1}z+\hbar^{-1}\frac{\theta_2^2}{2!}
p_{1,0}z+\hbar^{-1}\frac{\theta_2^3}{3!}(p_{2,1}+\hbar^{-1}\frac{1}{2!}p_{1,1}^2)z^2\\
&+&\frac{\theta_2^4}{4!}(2p_{2,0}+p_{1,1}p_{1,0})z^2+\dots
\end{eqnarray*}
where $\hbar^{-1}\theta_2 p_{1,1} z$ corresponds to a degree $1$
rational map with one interior marked point mapping to a specified
point in $\proj^2$ and one boundary marked point mapping to a
specified point in $L$; and
$\hbar^{-1}\frac{\theta_2^4}{4!}p_{1,1}p_{1,0}z^2$ corresponds to
a degree $2$ rational map through four specified, generic points
in $\proj^2$ with two boundary points of contact of order $1$ to
$L$, one at a specified point, the other free.
\end{example}

\subsection{The Rubber Potential}

We can write all intersection numbers arising from $\ma$ in terms
of a generating function.  Let $L$ be a line-bundle over a
projective manifold $X$.  Pick an Euler characteristic $\chi$, $m$
interior marked points, $r_0+r_\infty$ boundary marked points. Fix
a curve class $d\in B_1(X)$, a $r_0$-tuple of multiplicities
$(s^0_1,\dots,s^0_{r_0})$ to $D_0$ at the $r_0$ boundary marked
points, and a $r_\infty$-tuple of multiplicities
$(s^\infty_1,\dots,s^\infty_{r_\infty})$ to $D_\infty$ at the
$r_\infty$ boundary marked points. Let $\Xi$ be the set of rubber
graphs $\Gamma$ with no vertices associated to trivial cylinders
so that
\begin{enumerate}
\item[(1)] $\sum_v d(v)=d$

\item[(2)] $\sum_v (2-2g(v))=\chi$

\item[(3)] $|M|=m$

\item[(4)] $|R_0|=r_0$

\item[(5)]
$(\mu^0(1),\mu^0(2),\dots,\mu^0(r_0))=(s^0_1,s^0_2,\dots,s^0_{r_0})$

\item[(6)] $|R_\infty|=r_\infty$

\item[(7)]
$(\mu^\infty(1),\mu^\infty(2),\dots,\mu^\infty(r_\infty))=(s^\infty_1,s^\infty_2,\dots,s^\infty_{r_\infty}).$
\end{enumerate}
We have evaluation maps at the marked points
$$\Ev:\cm(\ca,\Gamma)\srarr X^m\times X^{r_0}\times X^{r_\infty}.$$
We also have the two line bundles on $\ma$, $\Top$ and $\Bot$.
Given a cohomology class $c\in H^*(X^{m+r_0+r_\infty})$, we
consider intersection numbers of the form
$$c_1(\Bot)^l \cup \Ev^*(c)\cap \vir{\ma}.$$

\begin{definition} Given a curve class $d\in B_1(X)$, an Euler
characteristic $\chi$, non-negative integers $l,n,r_0,r_\infty$,
multiplicities $S^0=(s^0_1,\dots,s^0_{r_0})$,
$S^\infty=(s^\infty_1,\dots,s^\infty_{r_\infty})$,
$$c_1,\dots,c_m,e^0_1,\dots,e^0_{r_0},e^\infty_1,\dots,e^\infty_{r_\infty}\in
H^*(X),$$ a non-negative integer $m$, define the {\em correlator}
\begin{eqnarray*}
&&(c_1,\dots,c_n \cdot e^0_1,\dots,e^0_{r_0} \cdot
e^\infty_1,\dots,e^\infty_{r_\infty})_{\chi,d,S^0,S^\infty,l}\\
&=&\sum_{\Gamma\in\Xi} \deg(c_1(\Bot)^l \cup
\Ev^*(c_1\times\dots\times c_n\times e^0_1\times\dots\times
e^0_{r_0}\times e^\infty_1\times \dots \times
e^\infty_{r_\infty})\\
&&\ \cap \vir{\cm(\ca,\Gamma)}.)
\end{eqnarray*}
\end{definition}

We organize the correlators into a generating function which we
call the {\em rubber potential} which takes values in a particular
graded algebra, $\ccr$. Pick a a homogeneous basis $c_1,\dots,c_k$
for $H^*(X)$. We define elements $z^A$ for all $A\in B_1(X)$,
$\hbar$, $\lambda$, $\beta_1,\dots,\beta_k$,
$p_{n,1},\dots,p_{n,k}$, $q_{n,1},\dots,q_{n,k}$ for all positive
integers $n$.  The $z^A$'s obey the relations $z^A \cdot z^B =
z^{A+B}$ where $+$ is addition in $B_1(X)$.  The elements are
graded of the following degrees
\begin{eqnarray*}
\deg z^d&=&2<c_1(TX),d>+2<c_1(L),d>\\
\deg \hbar&=&-2(\dim X-2)\\
\deg \lambda&=&-2\\
\deg \beta_i&=&2-\deg c_i\\
\deg p_{n,i}&=&2-\deg c_i-2n\\
\deg q_{n,i}&=&2-\deg c_i+2n.
\end{eqnarray*}

Multiplication in the algebra is defined as follows.  The
$\hbar$-, $\lambda$-, $\beta_i$-variables are taken to be
supercentral while the $p$- and $q-$variables obey
supercommutation relations
%%\begin{eqnarray*}
$$[p_{n_1,i_1},p_{n_2,i_2}]=0,\ [q_{n_1,i_1},q_{n_2,i_2}]=0,\
[q_{n_1,i_1},p_{n_2,i_2}]=n_1\delta_{n_1,n_2}g^{i_1 i_2}\hbar.$$
%%\end{eqnarray*}
where $g^{i_1 i_2}$ is the Poincare pairing on $H^*(X)$. Note that
this algebra can be realized by writing $q_{n,i}$ as a
differential operator
$$q_{n,i}=n\hbar\sum_j g^{ij}\frac{\partial}{\partial
p_{n,j}}.$$

\begin{figure}
\psfig{file=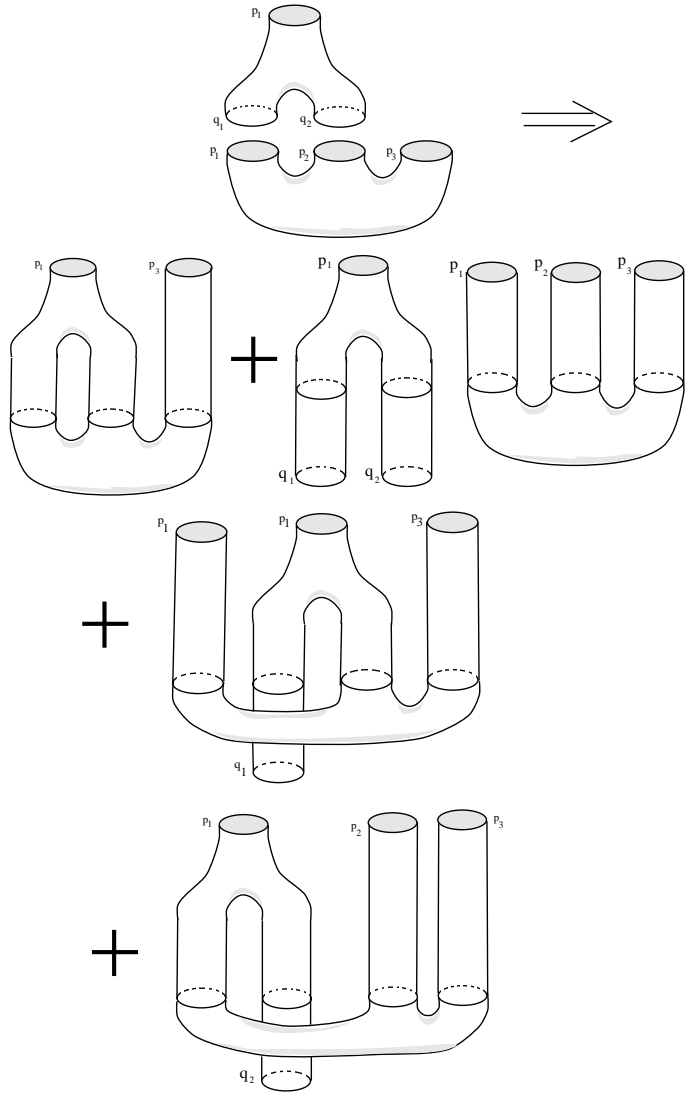,height=145mm} \caption{Geometric Illustration
of Multiplication in $\ccr$.  This figure is borrowed from
\cite{EGH}.} \label{figcomb}
\end{figure}

The multiplication keeps track of different ways of joining
curves. Let us consider an example in a toy model of our algebra.
Consider variables $p_1,p_2,p_3$ which are all of even parity
together with $q_1,q_2,q_3$ so
$$q_i=\hbar\frac{\partial}{\partial p_i}.$$
Then,
$$(\hbar^{-1}p_1q_1q_2)*(\hbar^{-1}p_1p_2p_3)=p_1p_3+\hbar^{-2}p_1^2p_2p_3q_1q_2
+\hbar^{-1}p_1^2p_3q_1+\hbar^{-1}p_1p_2p_3q_2.$$ If we see $\hbar$
as a genus marker where a term of Euler coefficient $\chi$ is
marked with $\hbar^{-\frac{1}{2}\chi}$, this multiplication
corresponds to the geometric situation illustrated in figure
\ref{figcomb}.

%%Note that if $g_{i_1 i_2}\neq 0$, $\deg c_{i_1}+\deg
%%c_{i_2}=2\dim_\com X$.  Therefore,
%%$$\deg(n_1\delta_{n_1,n_2}g^{i_1 i_2}\hbar)=-2(\dim X-2)=\deg
%%p_{n_1,i_1}+\deg q_{n_2,i_2}$$ and the bracket operation respects
%%the grading.

The algebra $\ccr$ consists of Laurent series in $\hbar$ whose
coefficients are power series in the $p$-variables whose
coefficients are power series in the $\beta$-variables whose
coefficients are polynomials in the $q$- and $\lambda$-variables.

\begin{definition} The {\em $\ccr$-correlator} is
\begin{eqnarray*}
&&(\beta_{i_1}c_{i_1}\dots\beta_{i_m}c_{i_m}\
,p_{m_1,j_1}c_{j_1}\dots p_{m_{r_0},j_{r_0}}c_{j_{r_0}},\
q_{n_1,k_1}c_{k_1}\dots
q_{n_{r_\infty},k_{r_\infty}}c_{k_{r_\infty}})_{\chi,d,l}\\
&=&\beta_{i_1}\dots\beta_{i_m}p_{m_1,j_1}\dots p_{m_{r_0},j_{r_0}}
q_{n_1,k_1}\dots q_{n_{r_\infty},k_{r_\infty}}\\
&&\cdot(c_{i_1},\dots,c_{i_n} \cdot c_{j_1},\dots,
c_{j_{r_0}}\cdot
c_{k_1},\dots,c_{k_{r_\infty}})_{\chi,d,(m_1,m_2,\dots,m_{r_0}),(n_1,n_2,\dots,n_{r_\infty}),l}
\end{eqnarray*}
considered as an element of $\ccr$.
\end{definition}

Let $\cb,\cq,\ccp$ be the (noncommutative) power-series algebras
freely generated by $\beta_i c_i$, $q_{n,i} c_i$, $p_{n,i}c_i$
respectively.

We can extend the $\ccr$-correlator by linearity to give a
multi-linear function
$$(\ ,\ ,\ )_{\chi,A,m}:\cb\otimes\ccp\otimes\cq\srarr\ccr$$

Let $B\in\cb$, $P\in\ccp$, $Q\in\cq$ be given by
\begin{eqnarray*}
B&=&\sum_{l\geq 0} \frac{1}{l!}\left(\sum_i \beta_i c_i\right)^l,\\
P&=&\sum_n \frac{1}{n!}\left(\sum_{k,i} p_{k,i}c_i\right)^n,\\
Q&=&\sum_n \frac{1}{n!}\left(\sum_{k,i} q_{k,i}c_i\right)^n.
\end{eqnarray*}

\begin{definition} \label{accpot} The {\em rubber potential} $A$
is
$$A=\sum_\chi \sum_{d\in B_1(X)} \sum_l \hbar^{-\frac{1}{2}\chi}\frac{\lambda^l}{l!} (B,P,Q)_{\chi,d,l} z^d.$$
\end{definition}

Note that for a moduli stack $\ma$ to be non-empty, by Lemma
\ref{multcond}, the multiplicities must satisfy
$$s_1^0+\dots+s_{r_0}^0-s_1^\infty-\dots-s_{r_\infty}^\infty=<c_1(L),A>$$
so
$$r_\infty\leq s_1^0+\dots+s^{r_0}_0-<c_1(L),A>.$$
Therefore, the rubber potential satisfies the polynomiality in $q$
condition to lie in $\ccr$. The rubber potential $A$ is
homogeneous of degree $2$.

\begin{definition} The {\em rubber potential without powers of
$\Bot$} is given by $A_{\lambda=0}$.
\end{definition}

\begin{example} The rubber potential of $(X,L)=(\pt,1_{\pt})$ obeys
$$A_{\lambda=0}=\hbar^{-1}(\frac{1}{3!}\theta_0^3+\frac{1}{2}\sum_{k,l\geq
1}(p_{k+l}q_k q_l+p_k p_l q_{k+l}))-\hbar^0\frac{1}{24}\theta_0.$$
Note that $\frac{1}{2}\sum_{k,l\geq 1}(p_{k+l}q_k q_l+p_k p_l
q_{k+l}))$ are the cut-and-join operators of \cite{GJV}.  The full
rubber potential can be related to Hurwitz numbers by use of a
localization argument in \cite{LZ}.
\end{example}

\begin{example} As a consequence of Corollary \ref{fourierrat},
for $(X,L)=(\proj^1,\oh(1))$, the terms in the rubber potential
without powers of $\Bot$ corresponding to positive degree maps are
\begin{eqnarray*}
A_{\lambda=0}&=&\hbar^{-1}\frac{1}{2\pi}\int_0^{2\pi}
\frac{(\beta_0+\sum_k p_{k,0}e^{-ikx}+\sum_k q_{k,0}e^{ikx})^2}{2}\\
&&\cdot (\beta_2+ \sum_k
p_{k,2}e^{-ikx}+\sum_k q_{k,2}e^{ikx}) dx\\
&+&\hbar^{-1}\frac{1}{2\pi}\int_0^{2\pi} e^{\beta_2+\sum_k
p_{k,2}e^{-ikx}+\sum_k q_{k,2}e^{ikx}} ze^{ix} dx.
\end{eqnarray*}
\end{example}

\subsection{Trivial Cylinders}

It was our convention to exclude trivial cylinders from the rubber
potential.  They will be accounted for by the algebra $\ccr$.  To
prove this, it will be advantageous to write down a potential
including trivial cylinders and relate it to the rubber potential.
Let $\Gamma$ be some rubber graph. Let $\Gamma_|$ be a rubber
graph obtained from $\Gamma$ by adjoining a degree $r$ trivial
cylinder.  From \cite{Ka2}, we have
\begin{theorem} There is a natural map
$$v:\cm(\ca,\Gamma_|)\srarr \cm(\ca,\Gamma)\times X$$
so that
$$v_*\vir{\cm(\ca,\Gamma_|)}=\frac{1}{r}\vir{\cm(\ca,\Gamma)}\times
[X]$$ and
$$v^*(\Bot)=\Bot.$$
{\rm Consequently if $\Gamma$ has $m$ interior marked points and
$r_0+r_\infty$ boundary marked points then we have a commutative
diagram
$$\xymatrix{
\cm(\ca,\Gamma_|)\ar[r]^>>>>>>>>>{\Ev_|}\ar[d]& X^m\times
X^{r_0+1}\times
X^{r_\infty+1}\ar[d]^h\\
\cm(\ca,\Gamma)\times X\ar[r]^<<<<<{\Ev\times\Delta}&X^m\times
X^{r_0}\times X^{r_\infty}\times (X\times X)}$$ where
$\Delta:X\srarr X^2$ is the diagonal and the morphism $h$ reorders
the products of $X$ so that the product of $X$'s corresponding to
the $r_0+1$st and $r_\infty+1$st boundary marked points are taken
to $X\times X$.}

For classes $A\in H^*(X^m\times X^{r_0+1}\times X^{r_\infty+1})$,
we have
\begin{eqnarray*}
& & \deg(\Ev_|^*(A)\cap\vir{\cm(\ca,\Gamma_|)})\\
&=&\frac{1}{r}\deg(((\Ev\times\Delta)\circ h^{-1})^*(A)\cap
(\vir{\cm(\ca,\Gamma)}\times [X])).
\end{eqnarray*}
\end{theorem}

\begin{definition} The {\em rubber potential with trivial
cylinders}, $A_|$ is defined as before except that we allow the
set $\Xi$ to contain graphs that have trivial cylinders for
vertices.
\end{definition}

Define the action of an algebra of power series in infinitely many
non-commuting variables $\kappa_1,\kappa_2,\dots$ on monomials
$f\in \ccr$ by
$$\kappa_n\cdot f=\frac{1}{n}\hbar^{-1}\sum_{i_1,i_2}(-1)^{(\deg_{pq}(f))(\deg(p_{n,i_1}))}g^{i_1 i_2}p_{n,i_1}f
q_{n,i_2}$$ where $\deg_{pq}(f)$ is the sum of the degrees of the
$p$ and $q$ variables in $f$. Extend the action linearly to
$\ccr$.  Let
$$TK=\sum_n \kappa_n$$
and define a map
$$\begin{array}{ccccc}
T & : & \ccr & \srarr &\ccr\\
T & : & f & \mapsto & e^{TK}f.
\end{array}$$

\begin{lemma} $T$ takes the rubber potential to the rubber
potential with trivial cylinders,
$$T(A)=A_|$$
\end{lemma}

\begin{proof} The proof is straightforward.  The factorial terms in the
exponential come from relabelling the boundary marked points.
\end{proof}

\begin{definition} Let $f$ and $h$ be elements in $\ccr$. We
define a binary operation $f *_| h$ as follows.  Introduce a set
of auxiliary variables $\tilde{p}_{n,i}$, $\tilde{q}_{n,i}=\sum_j
n \hbar g^{ij}\frac{\partial}{\partial \tilde{p}_{n,j}}$. Write
$f(p,\tilde{q})$, $h(\tilde{p},q)$, that is we substitute the
tilded variables into the power series. Define
$$f*_| h=f(p,\tilde{q})h(\tilde{p},q)|_{\tilde{p}=0}.$$
\end{definition}

Note that in the above we treat $\frac{\partial}{\partial
\tilde{p}_{n,j}}$ as an element with the same parity as $p_{n,j}$.
The operation, $*_|$ is the one that corresponds to stacking
curves to form multi-level curves.  This will be elaborated in the
section on degenerations.

\begin{lemma} $T$ is a homomorphism from $(\ccr,*)$ to
$(\ccr,*_|)$.\end{lemma}

\begin{proof}
$$\frac{\partial}{\partial p_{n,i}}(e^{TK}1)=\hbar^{-1}\frac{1}{n}\sum_j
(e^{TK}1)g^{ij}q_{n,j}$$ which implies for $f$, a monomial,
$$\frac{\partial}{\partial p_{n,i}}(e^{TK}f)=\hbar^{-1}\frac{1}{n}\sum_j
(-1)^{(\deg_{pq}f)(\deg
p_{n,i})|}(e^{TK}f)g^{ij}q_{n,j}+e^{TK}\frac{\partial f}{\partial
p_{n,i}}.$$

The lemma follows by induction on the number of $p$ and $q$
variables in $f$.
\end{proof}

%%\begin{proof} It suffices to show
%%$$T(fh)=T(f)*_|T(h)$$
%%for $f$ and $h$ monomials. In addition, it is clear that if the
%%conclusion holds for $f$ and $h$, then we have
%%\begin{eqnarray*}
%%T(p_{n,i}fh)&=&T(p_{n,i}f)*_|T(h)\\
%%T(fhq_{n,i})&=&T(f)*_|T(hq_{n,i})
%%\end{eqnarray*}
%%Therefore, we can assume that $f$ involves no $p$ variables and
%%$g$ involves no $q$ variables.
%%
%%Now, we induct on the sum of powers of $q$ in $f$ and powers of
%%$p$ in $g$:
%%\begin{eqnarray*}
%%T(fq_{n,i})*_|T(h)&=&\pm T(f)*_|(\sum_j n\hbar
%%g^{ij}\frac{\partial}{\partial
%%p_{n,j}}T(h))\\
%%&=&\pm T(f)*_|(T(\sum_j \hbar ng^{ij} \frac{\partial h}{\partial
%%p_{n,j}}))\\
%%&&\ \ \pm(-1)^{(\deg_{pq}h)(\deg q_{n,j})|}T(f)*_|(T(h)\sum_j
%%g^{ij}q_{n,j})\\
%%&=&T(f\sum_j \hbar ng^{ij} \frac{\partial h}{\partial
%%p_{n,j}})\\
%%&&\ \ +(-1)^{(\deg_{pq}h)(\deg q_{n,j})}T(fh)\sum_j g^{ij}q_{n,j}\\
%%&=& T(fq_{n,i}h)
%%\end{eqnarray*}
%%where $\pm$ denotes a possible sign for each monomial summand in
%%each term.
%%\end{proof}

\begin{definition} \label{ratpot} {The {\em rational potential} $A$
is
$$A_\text{rat}=\sum_{A\in B_1(X)} (\Gamma,P,Q)^\bullet_{g=0,A,m} z^A$$
where $(\ ,\ ,\ )^\bullet$ where the sum is taken over moduli
spaces $\ma$ involving only connected domains of genus 0.}
\end{definition}

\subsection{Action of $\ccr$ on $\cf$}

\begin{lemma} \label{gaction} {$\cf$ can be given the structure of
a graded $\ccr$-module.}\end{lemma}

\begin{proof}
Consider the inclusion
$$I:D\srarr Z$$
and the induced maps
$$\begin{array}{ccccc}
I^*&:&H^*(Z)&\srarr & H(D),\\
I^{*\vee}&:&H^*(D)^\vee & \srarr &H^*(Z)^\vee.
\end{array}$$

We define the action of $\ccr$ on $f\in\cf$ as follows:
\begin{eqnarray*}
\lambda\cdot f&=&0\\
\hbar\cdot f&=&\thbar f\\
p_{n,i}\cdot f&=&\tp_{n,i}\\
q_{n,i}\cdot f&=&\thbar n\sum_{i'} g^{ii'}\frac{\partial}{\partial
\tp_{n,i'}}f\\
\beta_i\cdot f&=&I^{*\vee}(\gamma_i)f\\
z^d\cdot f&=&\tz^{i_* d}f
\end{eqnarray*}
where $g^{ii'}$ is the intersection pairing on $H^*(D)$.

Because
$$\begin{array}{ccccc}
\deg \tz^{i_*
d}&=&<c_1(TZ),i_*d>&=&<i^*c_1(TZ),d>\\
&=&<c_1(TD)+c_1(N),d>&=&\deg z^d,
\end{array}$$
the action preserves grading.
\end{proof}

\begin{definition} Define a bilinear operation
$$\cdot_|:\ccr\otimes\cf\srarr \cf$$
as follows, for $f\in\cf$, $h\in\ccr$.  Given a monomial
$$h=\hbar^{-\frac{1}{2}\chi}\lambda^m\beta_{i_1}\dots\beta_{i_m}
p_{m_1,j_1}\dots p_{m_{r_0},j_{r_0}} q_{n_1,k_1}\dots
q_{n_{r_\infty},k_{r_\infty}}z^d\in\ccr$$ define $h\cdot_| f$ by
defining the action of $q_{n,i}$ on $\cf$ by
$$q_{n,i}=\left(\sum_{i'} g^{ii'}\frac{\partial}{\partial
\tilde{p}_{n,i'}}\right)$$ and defining
\begin{eqnarray*}
h\cdot_|f=\left((\thbar^{-\frac{1}{2}\chi}\delta_{m0}I^\vee(\beta_{i_1})\dots
I^\vee(\beta_{i_m}) p_{m_1,j_1}\dots p_{m_{r_0},j_{r_0}}
q_{n_1,k_1}\dots q_{n_{r_\infty},k_{r_\infty}})f\right)|_{\tp=0}
\end{eqnarray*}
and then by substituting $\tp_{n,j}$ for $p_{n,j}$
\end{definition}

This operation corresponds to joining a curve in $\ma$ to one in
$\mz$.  Analogously to the multiplication in $\ccr$, $\cdot_|$ and
the module structure $\cdot$ are related as follows:

\begin{lemma} \label{dotbar} For $h\in\ccr$, $f\in\cf$, we have
$$T(h)\cdot_|f=h\cdot f.$$
\end{lemma}

%\bibliographystyle{plain}
%\bibliography{thesis}

%\end{document}

%% file: dgr.tex
Theorem \ref{mpiobb} gives formulas relating the line-bundles
$$\Dil,\Split,\Lnt{i},\Lnb{i}$$
on $\cm(\cz,\Gamma_Z)$ and $\cm(\ca,\Gamma_A)$. In this section,
we will show that the first Chern classes of these line-bundles
turn out represent specific geometric situations involving split
curves. For example, $c_1(\Split)$ is a substack of
$\cm(\ca,\Gamma)$ that is, in a virtual sense, all split curves.
$c_1(\Lnt{i})$ virtually consists of all split curves in which the
$i$th marked point is not on the topmost component.  This allows
us to write the cap product of a first Chern class of one of our
bundles with the virtual cycle in terms of the virtual cycles of
smaller moduli spaces.  This provides {\em degeneration} formulae
that can be expressed in the language of generating functions.

We will express the first chern class of various line-bundles
geometrically by adapting Li's argument \cite{Li2}.  The argument
is in several stages and we state it only in the case
$\cm(\cz,\Gamma)$ noting that the case for $\cm(\ca,\Gamma)$ is
exactly analogous:
\begin{enumerate}
\item[(1)] For $\Gamma$, consider quadruples
$\Upsilon=(\Gamma_Z,\Gamma_A,L,J)$ so that the graph join,
$\Gamma_Z*_{L,J}\Gamma_A$ is isomorphic to $\Gamma$.  We can
define a line bundle $L_\Upsilon$ on $\cm(\cz,\Gamma_Z)$.

\item[(2)] We show that
$$c_1(L_\Upsilon)\cap \vir{\cm(\cz,\Gamma)}=m(\Upsilon)\vir{\cm(\cz\sqcup\ca,\Upsilon)}$$
where $\vir{\cm(\cz\sqcup\ca,\Upsilon)}$ is an appropriately
defined virtual cycle.

\item[(3)] Given the joining morphism
$$\Phi:\cm(\ca,\Gamma_A)\times_{D^r}
\cm(\cz,\Gamma_Z)\srarr\cm(\cz,\Gamma_A\sqcup_{L,J}\Gamma_Z)$$ and
the diagram
$$\xymatrix{\cm(\ca,\Gamma_A)\times_{D^r}
\cm(\cz,\Gamma_Z)\ar[r] \ar[d]
&\cm(\ca,\Gamma_A)\times \cm(\cz,\Gamma_Z)\ar[d]\\
D^r\ar[r]_\Delta & D^r\times D^r}$$ where $\Delta$ is the diagonal
map. We have
$$\Phi_*\Delta^!(\vir{\cm(\ca,\Gamma_A)}\times
\vir{\cm(\cz,\Gamma_Z)})=\vir{\cm(\ca\sqcup\cz,\Upsilon)}.$$

\item[(4)] Given a line-bundle $L=\Dil$ or $L=\Le{i}$ (or in the
case of $\ma$, $\Split$,$\Lnt{i}$,$\Lnb{i}$), we exhibit a set of
join-equivalence classes $\Omega$ so that
$$L=\otimes_{[\Upsilon]\in\Omega} L_\Upsilon.$$

\item[(5)] Consequently
$$c_1(L) \cap \vir{\cm(\cz,\Gamma)}=\sum_{\Upsilon\in\Omega}
m(\Upsilon)\Phi_*\Delta^!(\vir{\cm(\ca,\Gamma_A)}\times\vir{\cm(\cz,\Gamma_Z)}).$$
\end{enumerate}

To modify this argument to work for $\cm(\cz,\Gamma)$, replace all
pairs $(\Gamma_A,\Gamma_Z)$ with $(\Gamma_t,\Gamma_b)$ and replace
$\cz$ with $\ca$.

\subsection{Interpretation of Bundles}

Let us rewrite the bundles $\Dil$,$\Le{i}$,
$\Lnt{i}$,$\Lnb{i}$ as
tensor products of $L_\Upsilon$'s on $\cm(\cz,\Gamma)$ and
$\cm(\ca,\Gamma)$.

On $\mz$ where $i$ is the label for an interior marked point,
\begin{enumerate}
\item[(1)] $\Omega_{\Dil}=\{\Upsilon=(\Gamma_A,\Gamma_Z,L,J)\}$
the set of all join-equivalence classes of quadruples
$\Upsilon=(\Gamma_A,\Gamma_Z,L,J)$.

\item[(2)] $\Omega_{\Le{i}}=\{(\Gamma_A,\Gamma_Z,L,J)|i\in
J(M_A)\}.$
\end{enumerate}

while on $\ma$ where $i,j$ are labels for interior marked points,
\begin{enumerate}
\item[(1)] $\Omega_{\Split}=\{(\Gamma_t,\Gamma_b,L,J)\}.$

\item[(2)] $\Omega_{\Lnb{i}}=\{(\Gamma_t,\Gamma_b,L,J)|i\in
J(M_t)\}$.

\item[(3)] $\Omega_{\Lnt{i}}=\{(\Gamma_t,\Gamma_b,L,J)|i\in
J(M_b)\}$.

\item[(4)] $\Omega_{(i,j)}=\{(\Gamma_t,\Gamma_b,L,J)|i\in
J(M_t),j\in J(M_b)\}$.

\item[(5)] $\Omega_{(ij,)}=\{(\Gamma_t,\Gamma_b,L,J)|i,j\in
J(M_t)\}.$

\item[(6)] $\Omega_{(,ij)}=\{(\Gamma_t,\Gamma_b,L,J)|i,j\in
J(M_b)\}.$
\end{enumerate}

\begin{theorem} \label{interpofbundles} \cite{Ka2} {For
$L=\Dil,\Le{i},\Split,\Lnb{i},\Lnt{i}$,
$$L=\bigotimes_{[\Upsilon]\in \Omega_L} L_\Upsilon.$$
where $[\Upsilon]$ denotes a join-equivalence class and $\Upsilon$
a representative element.}
\end{theorem}

\subsection{Splitting of Moduli Stacks}

We need to cite a number of results from \cite{Li2}.  These
results were proved for a different moduli stack, $\cm(\cw)$, but
because of the explicit parallels between that space and the
construction of $\cm(\ca,\Gamma)$ and $\cm(\cz,\Gamma)$, the
proofs can be modified in straightforward fashion.  We begin by
relating the virtual cycle $\vir{\cm(\ca\sqcup\cz,\Upsilon)}$
defined in \cite{Li2} where $\Upsilon=(\Gamma_Z,\Gamma_A,L,J)$ is
a graph-join quadruple to other virtual cycles.

\begin{theorem} {We have the following equality among cycle classes
$$c_1(L_\Upsilon)\cap\vir{\cm(\cz,\Gamma_A*_{L,J}\Gamma_Z)}=m(\Upsilon)\vir{\cm(\ca\sqcup\cz,\Upsilon)}$$}
where $m(\Upsilon)$ is as in Definition \ref{multgraphjoin}.
\end{theorem}

Consider the fiber square
$$\xymatrix{\cm(\ca,\Gamma_A)\times_{D^r}
\cm(\cz,\Gamma_Z)\ar[r] \ar[d]
&\cm(\ca,\Gamma_A)\times \cm(\cz,\Gamma_Z)\ar[d]\\
D^r\ar[r]_\Delta & D^r\times D^r}$$ where $\Delta$ is the diagonal
morphism and the downward maps are induced from evaluation at the
boundary marked points of $\mz$ and the boundary marked points at
$D_\infty$ on $\ma$.  Let the virtual cycle on the fiber product
be given by
$$\vir{\cm(\ca,\Gamma_A)\times_{D^r}\cm(\cz,\Gamma_Z)}=\Delta^!(\vir{\cm(\ca,\Gamma_A)}\times\vir{\cm(\cz,\Gamma_Z)}).$$

\begin{theorem} {If
$$M_{[\Upsilon]}=\coprod_{(\Gamma_A',\Gamma_Z',L,J)\in[\Upsilon]}
\cm(\ca,\Gamma_A')\times_{D^r}\cm(\cz,\Gamma_Z')$$ is given the
virtual cycle of a disjoint union, then
$$\Phi_{[\Upsilon]}:M_{[\Upsilon]}\srarr \cm(\cz,\Gamma_A*_{L,J}\Gamma_Z)$$
gives
$$\Phi_{[\Upsilon]*}(\vir{M})=|MZ|!|MA|!(|RZ|!)^2\vir{\cm(\ca\sqcup\cz,\Upsilon)}$$}
\end{theorem}

Note that the multiplicity term is natural in light of Proposition
\ref{totalgluing}.

\begin{corollary}
$$c_1(L_\Upsilon)\cap\vir{\cm(\cz,\Gamma_A*_{L,J}\Gamma_Z)}=\frac{m(\Upsilon)}{|MZ|!|MA|!(|RZ!)^2}\Delta^!(M_\Upsilon).$$
\end{corollary}

$L$ together with $i:X\srarr Z$ induces a morphism
$$\Lambda:(X^{|MA|}\times X^{|RA_0|}) \times Z^{|MZ|}\srarr
Z^M\times X^R$$ where $M=|MZ|+|MA|$ and $R=|RA_0|$ are the number
of interior and boundary marked points in $\Gamma_A *_{L,J}
\Gamma_Z$

We have morphisms
$$\xymatrix{
X^{|MA|+|RA_0|}\times X^{|RZ|} \times Z^{|MZ|}
\ar[r]^>>>>{\tilde{\Delta}} \ar[d]^p & X^{|MA|+|RA_0|} \times
X^{|RA_\infty|}
\times Z^{|MZ|} \times X^{|RZ|}\\
(X^{|MA|}\times X^{|RA_0|}) \times Z^{|MZ|}&& }$$ where
$\tilde{\Delta}$ is induced by $\Delta:X^{|RZ|}\srarr
X^{|RA_\infty|}\times X^{|RZ|}$ and $p$ is the projection.

Therefore, for $c\in H^*(Z^{|MZ|}\times X^{|RZ|}$,
$$\deg((\Ev^*(c)\cup
c_1(L_\Upsilon))\cap\vir{\cm(\cz,\Gamma_A*_{L,J}\Gamma_Z)})$$
$$=\frac{m(\Upsilon)}{\Aut_{\Gamma_Z,\Gamma_A,L}(RZ,RA_\infty)}
\deg(\Ev^*({\tilde{\Delta}}_!(p^*\Lambda^*c)\cap(\vir{\cm(\ca,\Gamma_A)}\times\vir{\cm(\cz,\Gamma_Z)})$$

%%We can abuse notation and write
%%$$\deg((\Ev^*(C)\cup
%%c_1(L_\Upsilon))\cap\vir{\cm(\cz,\Gamma_Z*_{L,J}\Gamma_A)})$$
%%$$=\frac{m(\Upsilon)}{|MZ|!|MA|!|RZ!|^2}\deg(\Ev^*({\tilde{\Delta}}_!(p^*\Lambda^*C))
%%\cap \vir{P_\Upsilon})$$

%%Again, there are obvious degeneration formulae on
%%$\cm(\ca,\Gamma_t*_{L,J}\Gamma_b)$ obtained by replacing
%%$\Gamma_Z$ with $\Gamma_b$ and $\Gamma_A$ with $\Gamma_t$.

If $\Omega$ is one of the sets of join-equivalence classes, from
$$c_1(L_\Omega)=\sum_{[\Upsilon]\in\Omega} c_1(L_\Upsilon)$$
we have
\begin{theorem} \label{degenformula}
$$
\deg((\Ev^*(c)\cup c_1(L_\Omega))\cap\vir{\cm(\cz,\Gamma)}) =$$
$$
\sum_{[\Upsilon]\in \Omega}
\frac{m(\Upsilon)}{|MZ|!|MA|!|RZ|!^2}\deg(\Ev^*({\tilde{\Delta}}_!(p^*\Lambda^*c))
\cap \vir{M_\Upsilon})$$ and analogously for $\ma$.
%%where $c\in
%%H^*(X^m\times X^{r_0} \times X^{r_\infty})$,
%%$$\deg((\Ev^*(c)\cup
%%c_1(L_\Omega))\cap\vir{\cm(\ca,\Gamma)}) =$$
%%$$\sum_{[\Upsilon]\in \Omega}
%%\frac{m(\Upsilon)}{|MA_b|!|MA_t|!|{RA_{b0}}|!^2}\deg(\Ev^*({\tilde{\Delta}}_!(p^*\Lambda^*c))
%%\cap \vir{M_\Upsilon})
%%$$
\end{theorem}

\begin{definition} Let the
symbols \noindent
\begin{center}
\setlength{\unitlength}{.5cm}
\begin{picture}(12,5)
  \put(2,3.15){\line(1,-3){.5}}
  \put(2,2.85){\line(1, 3){.5}}
  \put(2.1,2.2) {\line(1,0) {.5}}
  \put(2.1,3.8) {\line(1,0) {.5}}
  \put(1.7,2)   {$j$}
  \put(1.7,3.6)   {$i$}

  \put(3.3,1.6) {$,$}
  \put(6,2.85)    {\line(1, 3){.5}}
  \put(6,3.15)    {\line(1, -3){.5}}
  \put(6.05,3.5)  {\line(1,0) {.5}}
  \put(6.1,3.95)  {\line(1,0) {.5}}
  \put(6.6,3.2)   {$j$}
  \put(5.65,3.75) {$i$}

  \put(7.1,1.6) {$,$}
  \put(10,2.85)     {\line(1, 3){.5}}
  \put(10,3.15)   {\line(1,-3){.5}}
  \put(10.05,2.05){\line(1,0) {.5}}
  \put(9.95,2.5)  {\line(1,0) {.5}}
  \put(10.6,1.8)  {$j$}
  \put(9.55,2.25){$i$}
\end{picture}
\end{center}
denote the cohomology classes
$$\sum_{[\Upsilon]\in\Omega_{(i,j)}} c_1(L_\Upsilon),\sum_{[\Upsilon]\in\Omega_{(ij,)}}
c_1(L_\Upsilon),\sum_{[\Upsilon]\in\Omega_{(,ij)}}
c_1(L_\Upsilon)$$ respectively.
\end{definition}

These cohomology classes are dual to the cycles in $\ma$
representing split curves with $i$ and $j$ specified on top and
bottom component as specified in the symbol, counted with the
appropriate weight.

\subsection{Normal Bundle to Split Curves}

The following is useful for localization computations.

Let $(\Gamma_Z,\Gamma_A,L,J)$ be graph-join quadruple.  Consider
$L_\Upsilon$ for the quadruple $\Upsilon=(\Gamma_Z,\Gamma_A,L,J)$.
Then $c_1(L_\Upsilon)$ is a substack of
$\cm(\cz,\Gamma_Z*_{L,J}\Gamma_A)$.

Consider the moduli stacks $\cm(\cz,\Gamma_Z)$,
$\cm(\ca,\Gamma_A)$, $\cm(\cz,\Gamma_A*_{L,J}\Gamma_Z)$, and  the
inclusion
$$\Phi:\cm(\ca,\Gamma_A)\times_{D^r} \cm(\cz,\Gamma_Z)\srarr \cm(\cz,\Gamma_A*_{L,J}\Gamma_Z).$$
$\cm(\ca,\Gamma_Z)\times_{D^r} \cm(\cz,\Gamma_A)$ has projections
$p_A,p_Z$ to its $\ma$ and $\mz$ factors.

\begin{theorem} On $\cm(\cz,\Gamma_A*_{L,J}\Gamma_Z)$, {$\Phi^*L_\Upsilon=p_Z^*\Dil^\vee \otimes
p_A^*\Bot^\vee$.}\end{theorem}

Similarly, we have
\begin{theorem} On $\cm(\ca,\Gamma_{A_t}*_{L,J}\Gamma_{A_b})$, {$\Phi^*L_\Upsilon=p_{Ab}^*\Top^\vee \otimes
p_{At}^*\Bot^\vee$.}\end{theorem}

\subsection{Degeneration Formulae}

The above degeneration formulas can be written in terms of
generating functions.

For $L=\Dil$, we can write down a potential $F_{\Dil}$ which is
defined by a formula similar to that of the relative potential
except that instead of evaluating all possible cohomology classes
on $\vir{\mz}$, we evaluate them on $c_1(\Dil)\cap \vir{\mz}$.
That is, we define the $\Dil$ correlator by for $a_1,\dots,a_n\in
H^*(Z)$, $b_1,\dots,b_r\in H^*(X)$
$$<a_1,\dots,a_m\cdot b_1,\dots,b_r>_{\Dil,\chi,A,(s_1,\dots,s_r)}$$
$$=\sum_{\Gamma\in\Xi} (\Ev^*(a_1\times\dots\times a_m\times
b_1\times\dots\times b_r)\cup c_1(\Dil))\cap
\vir{\cm(\cz,\Gamma)}$$ and define $F_\Dil$ as in the previous
section. Define the rubber potential with trivial cylinders
without powers of $c_1(\Bot)$ by
$${A_|}_{\lambda=0}=(A_|)_{\lambda=0}.$$
Then we have by Lemma \ref{dotbar}, Theorem \ref{interpofbundles},
and Theorem \ref{degenformula}
\begin{theorem}
\begin{eqnarray*}
F_{\Dil}&=&(A_|)_{\lambda=0} \cdot_| F\\
F_{\Dil}&=&A_{\lambda=0} \cdot F.
\end{eqnarray*}
\end{theorem}

To study insertions of $c_1(\Le{i})$, we choose an element $e_j$
of our basis for $H^*(Z)$.  Because $c_1(\Le{i})$ is dependent on
the choice of marked point, we add a {\em distinguished} marked
point to all of the relative graphs that contribute to our
potential.  At this marked point, we evaluate
$c_1(\Le{i})\cup\ev_i^*e_j$.  More formally, given a graph
$\Gamma$ with $m$ marked points, consider the set $D(\Gamma)$
consisting of all graphs $\Gamma'$ with $m+1$ marked points such
that when we forget $m+1$st marked point on $\Gamma'$, we obtain
$\Gamma$.  Consider the $(\text{ex},e_j)$ correlator given by
$a_1,\dots,a_m\in H^*(Z)$, $b_1,\dots,b_r\in H^*(X=D)$
\begin{eqnarray*}
&<a_1,\dots,a_m\cdot
c_1,\dots,c_r>_{(\text{ex},e_j),\chi,A,(s_1,\dots,s_r)}\\
=&\sum_{(\Gamma,k)\in\Xi}(\sum_{\Gamma'\in D(\Gamma)}
(\Ev^*(a_1\times\dots\dots\times a_m\times e_j\times
b_1\times\dots\times b_r)\\
&\cup c_1(\Le{m+1})\cap \vir{\cm(\cz,\Gamma')}))
\end{eqnarray*}
We write down $F_{\text{ex},j}$ by using the modified correlator.
Then, for $i:X=D\srarr Z$, we write
$$i^*e_j=\sum_{l}  M_{jl} c_l$$
for $M_{jl}\in\rat$ which gives
\begin{eqnarray*}
F_{\text{ex},j}&=&\sum_l M_{jl} \frac{\partial
(A_|)_{\lambda=0}}{\partial \beta_l}\cdot_| F\\
&=&\sum_l M_{jl} \frac{\partial A_{\lambda=0}}{\partial
\beta_l}\cdot F.
\end{eqnarray*}
Now, since $\Le{k}=\ev_k^*([D])$, we have

\begin{theorem} \label{pointsfromD} Let $e_j$ be an element of our
basis for $H^*(Z)$.  Let $N$ be a matrix defined by
$$e_j\cup [D]=\sum_l N_{jl} e_l$$
%%$$i^*e_j=\sum_{l} M_{jl} c_l$$
then
$$\sum_l N_{jl}\frac{\partial F}{\partial
\theta_l}=\sum_l M_{jl} \frac{\partial A_{\lambda=0}}{\partial
\beta_l}\cdot F$$
\end{theorem}

\begin{proof} $\ev_k^*e_j\cup c_1(\Le{k})=(\ev_k^*(e_j\cup
[D]))$ implies
$$F_{\text{ex},j}=\sum_l N_{jl}\frac{\partial F}{\partial
\theta_l}.$$
\end{proof}

Likewise, we can write down a rubber potential with $c_1(\Split)$
inserted in the correlator.  This is analogous to $\Dil$ in the
relative case.

\begin{theorem} {$A_{\Split}=A_{\lambda=0} * A.$}
\end{theorem}

We can write down a potential involving insertions of
$c_1(\Lnb{i})$.  This is analogous to $\Le{i}$ in the relative
case.  Again, we have to single out a cohomology class $c_j\in
H^*(X)$ where at some marked point $i$, we will evaluate
$\ev_i^*(c_j)$ and $\Lnb{i}$.

\begin{theorem}{$A_{\Lnb{i},c_j}=\frac{\partial
A_{\lambda=0}}{\partial \beta_j} * A.$}
\end{theorem}

Similarly, for $\Lnt{i}$, we have

\begin{theorem}\label{linottop}{$A_{\Lnt{i},c_j}=A_{\lambda=0}*\frac{\partial
A}{\partial \beta_j}.$}
\end{theorem}

From Theorem \ref{mpiobb}(3) one can obtain a degeneration formula
for the rubber potential by inserting
$$c_1(\Bot\otimes\ev_i^*L)\ev_i^*(c_j)=(c_1(\Bot)+ev_i^*c_1(L))\ev_i^*(c_j)$$
at a distinguished point in the rubber potential to compute
$A_{\Lnt{i},c_j}$.  Note that $\frac{\partial A}{\partial
\lambda}$ is the rubber potential with an extra insertion of
$c_1(\Bot)$.

\begin{theorem} \label{degenma} {Define the matrix $N_{ij}$ by
$$c_1(L)\cup c_i=\sum_j N_{ij} c_j.$$
Then for each $i$, we have
$$\frac{\partial}{\partial \lambda}\frac{\partial A}{\partial \beta_i}+\sum_j N_{ij}\frac{\partial
A}{\partial \beta_j}=\frac{\partial A_{\lambda=0}}{\partial
\beta_i}*A.$$}
\end{theorem}

\subsection{Reconstruction of the Relative Potential}

Let us look at the relative case with $(Z,D)$.  Let $N$ be the
normal bundle to $D$ in $Z$. We can use the degeneration formula
for $\Le{i}$ to reconstruct the relative Gromov-Witten invariants
of $(Z,D)$ from the rubber potential (without powers of $\Bot$) of
$(D,N)$ and from seed values of the relative Gromov-Witten
invariants of $(Z,D)$.  This is a formal consequence of the fact
that $\ev_i^*D=\Le{i}$ and Theorem \ref{pointsfromD}.

We need to pick a particular basis for $H^*(Z)$.  Let $V\subseteq
H^*(Z)$ be the subspace
$$V=\im(\cup [D]:H^{*-2}(Z)\srarr H^*(Z)).$$
Let $e_1,\dots,e_v$ be a homogeneous basis for $V$, ordered by
degree. Extend this to a basis $\{e_{v+1},\dots,e_{v+w}\}$ of
$H^*(Z)$.

\begin{theorem} \label{Freconstruct} The relative potential $F$ of
$(Z,D)$ can be reconstructed from the rubber potential of $(D,N)$
together with the relative potential involving only the classes
$\{e_{v+1},\dots,e_{v+w}\}$, that is, from
$$F|_{\theta_1=\theta_2=\dots=\theta_v=0}.$$
\end{theorem}

\begin{proof}  We add in one $\theta_i$ variable at a time.  So, suppose we
have determined
$$F|_{\theta_{j}=\dots=\theta_v=0}.$$
Since $e_j\in V$,
$$e_j=\sum_k a_k e_k\cup [D]$$
for $a_k\in\rat$ where $a_l=0$ for $l\in\{j,\dots,v\}$ for degree
reasons.
%%Define $M$, $N$ as above, by
%%\begin{eqnarray*}
%%i^*e_k&=&\sum_{l} M_{kl} c_l\\
%%e_k\cup [D]&=&\sum_l N_{kl} e_l
%%\end{eqnarray*}
%%Therefore,
%%$$\sum_k a_lN_{km}=\delta_{jm}$$
%%Therefore,
%%$$\sum_k a_k \sum_l N_{kl}\frac{\partial F}{\partial
%%\theta_l}=\sum_k a_k \sum_l M_{kl} \frac{\partial
%%A_{\lambda=0}}{\partial \beta_l}\cdot F$$ implies
By the above, we have
$$\frac{\partial F}{\partial \theta_j}=\sum_{k,l} a_kM_{kl}\frac{\partial
A_{\lambda=0}}{\partial \beta_l}\cdot F$$ which allows us to solve
for
$$F|_{\theta_{j+1}=\dots=\theta_v=0}.$$
\end{proof}

\subsection{Transferring Classes between Split Curves}
\begin{theorem} \label{tcbsc} On $\ma$ with two distinguished
interior marked points, the following equation holds among
divisors.

\noindent
\begin{center}
\setlength{\unitlength}{.5cm}
\begin{picture}(17,5)
  \put(1,3)  {$\ev_2^*(c_1(L))-\ev_1^*(c_1(L))$}
  \put(9.5,3)  {$=$}
  \put(12,3.15){\line(1,-3){.5}}
  \put(12,2.85){\line(1, 3){.5}}
  \put(12.1,2.2) {\line(1,0) {.5}}
  \put(12.1,3.8) {\line(1,0) {.5}}
  \put(11.7,2)   {$2$}
  \put(11.7,3.6)   {$1$}
  \put(14,3)   {$-$}
  \put(16,3.15){\line(1,-3){.5}}
  \put(16,2.85){\line(1, 3){.5}}
  \put(16.1,2.2) {\line(1,0) {.5}}
  \put(16.1,3.8) {\line(1,0) {.5}}
  \put(15.7,2)   {$1$}
  \put(15.7,3.6)   {$2$}
\end{picture}
\end{center}
\end{theorem}

\begin{proof} Recall the following facts:
\begin{eqnarray*}
c_1(\Top)+c_1(\Bot)&=&c_1(\Split),\\
c_1(\Top)&=&c_1(\Lnt{i})+\ev_i^*c_1(L),\\
c_1(\Bot)&=&c_1(\Lnb{i})-\ev_i^*c_1(L).;
\end{eqnarray*}

Putting these together, we see
$$c_1(\Split)=c_1(\Lnt{1})-ev_1^*c(L)+c_1(\Lnt{2})+ev_2^*c_1(L).$$
Diagrammatically, this equation can be expressed as

\noindent \setlength{\unitlength}{.45cm}
\begin{center}
\begin{picture}(28,11)
  \put(1,7.85)    {\line(1, 3){.5}}
  \put(1,8.15)    {\line(1, -3){.5}}
  \put(1.05,8.5)  {\line(1,0) {.5}}
  \put(1.1,8.95)  {\line(1,0) {.5}}
  \put(1.6,8.2)   {$2$}
  \put(0.65,8.75) {$1$}
  \put(2.5,8)     {$+$}
  \put(4,8.15)    {\line(1,-3){.5}}
  \put(4,7.85)    {\line(1, 3){.5}}
  \put(4.1,7.2)   {\line(1,0) {.5}}
  \put(4.1,8.8)   {\line(1,0) {.5}}
  \put(3.7,7)     {$2$}
  \put(3.7,8.6)   {$1$}
  \put(5.5,8)     {$+$}
  \put(7,8.15)    {\line(1,-3){.5}}
  \put(7,7.85)    {\line(1, 3){.5}}
  \put(7.1,7.2)   {\line(1,0) {.5}}
  \put(7.1,8.8)   {\line(1,0) {.5}}
  \put(6.7,7)     {$1$}
  \put(6.7,8.6)   {$2$}
  \put(8.6,8)     {$+$}
  \put(10,7.85)     {\line(1, 3){.5}}
  \put(10,8.15)   {\line(1,-3){.5}}
  \put(10.05,7.05){\line(1,0) {.5}}
  \put(9.95,7.5)  {\line(1,0) {.5}}
  \put(10.6,6.8)  {$2$}
  \put(9.55,7.25){$1$}

  \put(1.5,3)    {$=$}

  \put(3,3.15)    {\line(1,-3){.5}}
  \put(3,2.85)    {\line(1, 3){.5}}
  \put(3.1,2.2)   {\line(1,0) {.5}}
  \put(3.1,3.8)   {\line(1,0) {.5}}
  \put(2.7,2)     {$1$}
  \put(2.7,3.6)   {$2$}
  \put(4.5,3)     {$+$}

  \put(6,2.85)     {\line(1, 3){.5}}
  \put(6,3.15)   {\line(1,-3){.5}}
  \put(6.05,2.05){\line(1,0) {.5}}
  \put(5.95,2.5)  {\line(1,0) {.5}}
  \put(6.6,1.8)  {$2$}
  \put(5.55,2.25){$1$}

  \put(7,3)      {$-\ ev_1^*c_1(L)\ +$}

  \put(12.5,3.15)  {\line(1,-3){.5}}
  \put(12.5,2.85)  {\line(1, 3){.5}}
  \put(12.5,2.2)   {\line(1,0) {.5}}
  \put(12.6,3.8)   {\line(1,0) {.5}}
  \put(12.2,2)     {$1$}
  \put(12.2,3.6)   {$2$}
  \put(13.5,3)     {$+$}

  \put(15.5,2.85)  {\line(1, 3){.5}}
  \put(15.5,3.15)  {\line(1, -3){.5}}
  \put(15.55,3.5)  {\line(1,0) {.5}}
  \put(15.6,3.95)  {\line(1,0) {.5}}
  \put(16.1,3.2)   {$2$}
  \put(15.05,3.75) {$1$}
  \put(16.7,3)      {$+\ ev_2^*c_1(L)$}
\end{picture}
\end{center}
\noindent The result follows. \end{proof}

One can use this result to reconstruct all rubber invariants from
those where all powers of $c_1(L)$ are at a single interior marked
point.

%\bibliographystyle{plain}
%\bibliography{thesis}

%\end{document}

%% file: examples.tex
\subsection{Computation of the rational rubber potential}

We compute the rational rubber potential of $L=\oh(m)$ over
$\proj^r$. We do this by comparing the rubber potential to the
rational Gromov-Witten invariants of $\proj^r$.  There is a clear
geometric reason why this should be possible.  Given a smooth
rational curve $\proj^1$ with interior and boundary marked points,
$\{x_1,\dots,x_k,p^0_1,\dots,p^0_{s^0},p^\infty_1,\dots,p^\infty_{s^\infty}\}$
with the multiplicities
$(m^0_1,\dots,m^0_{s^0},m^\infty_1,\dots,m^\infty_{s^\infty})$,
specifying a rubber map to $L$ is equivalent to finding a degree
$d$ map of $f:\proj^1\srarr \proj^r$  and a nowhere zero section
(defined up to $\com^*$-action) of
$$f^*L\otimes\oh(-m^0_1
p^0_1-\dots-m^0_{s^0}p^0_{s^0}+m^\infty_1 p^\infty_1+\dots
+m^\infty_{s^\infty}).$$ The numerical condition in Lemma
\ref{multcond} for multiplicities implies that this bundle must be
trivial. Therefore, there is only one section up to multiplication
by an element of $\com^*$. The automorphism group of the map to
rubber is equal to that of the stable map, so the rubber invariant
should equal the Gromov-Witten invariant.

Unfortunately, this intuitive picture may not be true for singular
curves. One, however is able to prove these results when the
target is $\proj^r$ in which case rubber invariants count maps of
smooth curves.

Let $L=\oh(m)$ be a line-bundle over $\proj^r$.  Let
$P=\proj_{\proj^r}(L\oplus 1_{\proj^r})$ be the projective
completion of $L$.  We consider a stack $\ma$ of rubber maps to
$(X,L)$.

To prove that intersections on $\ma$ occur away from singular
curves, we will use a Kleiman-Bertini theorem argument.  Consider
the commutative diagram

$$\xymatrix{
\ma \ar@{->}[drr]^{\Ev_{\ma}} \ar[d]^\ft &&&\\
\cmbar_{0,n+r_0+r_\infty}(\proj^r,d)\ar[rr]_>>>>>>>>>{\Ev_{\cmbar}}&&(\proj^r)^{n+r_0+r_\infty}}$$
%\noindent
%\setlength{\unitlength}{.5cm}
%\begin{center}
%\begin{picture}(22,8)
%  \put(3,6){$\ma$}
%  \put(3.9,5.5){\vector(0,-1){2.4}}
%  \put(2.7,4){\small{ft}}
%  \put(1,2){$\cmbar_{0,n+r_0+r_\infty}(\proj^r,d)$}
%  \put(8,2.2){\vector(1,0){7.4}}
%  \put(10,1){\small $\Ev_{\cmbar}$}
%  \put(6,6.2){\vector(3,-1){9.5}}
%  \put(10.2,5.2){\small $\Ev_{\ma}$}
%  \put(16,2){$(\proj^r)^{n+r_0+r_\infty}$}
%%  \put(3,3){$\mm\prt$}
%%  \put(3,1){$\mm$}
%%  \put(1,5){$L$}
%%  \put(1.1,4.9){\vector(0,-1){1.6}}
%%  \put(1.2,4.9){\vector(1,-1){1.75}}
%%  \put(3.2,2){$\pi$}
%%  \put(1.55,3.15){$i_F$}
%%  \put(.9,4){$I$}
%%  \put(2.1,4.1){$i_L$}
%%  \put(1.9,2.7){$\pi_L$}
%\end{picture}
%\end{center}
where $\ft$ is the map that takes a rubber map to
$\proj_{\proj^r}(\oh(n)\oplus \oh)$ to its projection to
$\proj^r$, contracts unstable components, and sees boundary and
interior marked points as marked points.

%%Given a cycle $K\subset X^n\times X^{r_\infty} \times X^{r_0}$ of
%%codimension equal to $\dim \ma$, we can perturb $K$ and show by a
%%dimension count that $\Ev_{\cmbar}^{-1}(K)$ is disjoint from the
%%singular maps in $\cmbar(\proj^r,d)$.  Then $\Ev_\ma^{-1}(k)$ is
%%disjoint from singular maps in $\ma$.

We need the following well-known:

\begin{theorem} {$\mrn(\proj^r,d)$ is non-singular and
$\ev_i:\mrn(\proj^r,d)\srarr \proj^r$ is a smooth morphism.}
\end{theorem}

\begin{theorem} The locus of maps with singular domain in
$\mrn(\proj^r,d)$ is the finite union of sub-stacks $M'$ so that
\begin{enumerate}
\item $M'$ is non-singular and of positive codimension.

\item The evaluation map $\Ev:M'\srarr X^{n+r_0+r_\infty}$ is
smooth.
\end{enumerate}
\end{theorem}

Now, we will use the Kleiman-Bertini theorem in the following
form.

%\tpoint{Theorem} {Let $f:X\srarr (\proj^r)^l$ be a morphism from
%$X$, a non-singular stack  of dimension $d$.  Then there is an
%open set $U\subseteq ((\proj^r)^\vee)^l$ so that if
%$H_1\times\dots\times H_l$ is the product of hyper-planes
%corresponding to a point in $U$, then the inverse image
%$f^{-1}(H_1\times\dots\times H_l)$ is either empty or non-singular
%of dimension $d-l$.

\begin{theorem} {Let $f:X\srarr (\proj^r)^p$ be a morphism from
$X$, a non-singular stack of dimension $d$.  Let
$c=[H]^{a_1}\times \dots \times [H]^{a_p}\in H^l((\proj^r)^p)$
where $[H]$ is a hyperplane class. Then there exists an open set
$U\subseteq ((\proj^r)^\vee)^l$ each point of which corresponds to
a product of linear subspaces $K=V_1\times \dots \times V_p$,
Poincare-dual to $c$, such that $f^{-1}(K)$ is either empty or
non-singular of dimension $\dim X-l$.}
\end{theorem}

%Let $\Delta_k:\proj^r\srarr (\proj^r)^k$ be the multi-diagonal
%embedding. Consider the map
%$$\Delta=\Delta_{a_1}\times\Delta_{a_2}\times\dots\times\Delta_{a_p}:(\proj^r)^p\srarr
%(\proj^r)^{a_1}\times\dots\times(\proj^r)^{a_p}=(\proj^r)^l.$$
%$c=\Delta^*(H^l)$.  Consider the composition $\Delta\circ
%f:X\srarr (\proj^r)^l$. Then Kleiman-Bertini provides a product of
%hyper-planes $H=H_1\times \dots \times H_l$ such that
%$(\Delta\circ f)^{-1}(H)$ is empty or non-singular of dimension
%$d-l$.  By shrinking $U$, we can suppose $\Delta$ is transversal
%along $H$. Therefore, $C=\Delta^{-1}(H)$ represents the
%Poincare-dual to $c$.

%\end{proof}

%%Let us now compute the rubber invariants involving possibly
%%disconnected curves, none of whose components are trivial
%%cylinders.

Let us first show that there is no contribution coming from
rational curves with disconnected domain.

\begin{theorem} \label{disconnect} Let $\ma$ be some rubber moduli
space consisting of maps from the disjoint union of $q\geq 2$
rational curves, none of which are trivial cylinders. Let $l=\dim
\ma$. If $c\in H^l((\proj^r)^{n+r_0+r_\infty})$, then there exists
$K$, a product of linear subspaces in
$(\proj^r)^{n+r_0+r_\infty}$, Poincare-dual to $c$ so that
$\Ev_\ma^{-1}(K)$ is empty and therefore,
$$\deg(\Ev_\ma^*(c)\cap \vir{\ma})=0.$$
\end{theorem}

\begin{proof}

Let us write
$\cmbar=\cmbar^\circ_{\chi,n+r_\infty+r_0}(\proj^r,d)$ for the
moduli space of stable maps of the disjoint union of $q$ rational
curves to $\proj^r$ of degree $d$.  $\chi=2q$ is the Euler
characteristic of the domain. Then, by the dimension formula,
$$\dim \ma-\dim \cmbar=\frac{1}{2}\chi-1=q-1>0.$$
Therefore, the codimension of the class $c\in
H^l((\proj^r)^{n+r_0+r_\infty})$ exceeds the dimension of
$\cmbar$, so its Poincare-dual can be represented by a product of
hyper-planes, $C$ so that that $\Ev_\cmbar^{-1}(C)$ is empty.
\end{proof}

Now, let us compute the rational rubber invariants coming from
curves with connected domains.  We begin with the straight-forward
combinatorial lemma

\begin{lemma} {A map in $\ma$ with connected rational singular
domain is mapped by $\ft$ to a map with singular domain in
$\cmbar(\proj^r,d)$}
\end{lemma}

%%\begin{proof} Suppose we have a map from a singular curve, $f:C\srarr\ _n
%%P$. Consider the composition $(\pi\circ f):C\srarr\ _n P\srarr
%%\proj^r$ where $P=\proj_\proj^r(L\oplus 1)$.  We have to show that
%%once we contract the components of $C$ that are unstable under
%%\pi\circ f$ then the curve is still singular.  An unstable
%%irreducible component must be mapped to a point under $\pi\circ f$
%%and have fewer than $3$ special points. Therefore, it must lie in
%%a fiber of $P\srarr \proj^r$ with at most two special points.
%%These special points must correspond to boundary marked points
%%with contact to to $D_0$ and $D_\infty$. Otherwise, since the sum
%%of multiplicities at $D_0$ equals to the sum of multiplicities to
%%$D_\infty$, the curve must be disjoint from $D_0$ and $D_\infty$,
%%so must be mapped to a point, and thus unstable as a map in $\ma$.
%%Therefore, there must be one point of tangency to $D_0$ and one
%%point of tangency to $D_\infty$. Therefore the component must map
%%to a fiber and be totally ramified over $D_0$ and $D_\infty$. Such
%%a component is called a {\em trivial cover}.

%%We look at two cases.  If $k\geq 1$, then there is a component
%%that is not a trivial cover that is mapped to each copy of $P$ in
%%$\ _n P$. So there are at least two components that are stable as
%%maps to $X$.

%%If $k=1$, there cannot be any trivial covers in curves with
%%connected domains. \end{proof}

Putting everything together,

\begin{theorem} {Let $\ma$ be some rubber moduli space with target
$(X,L)$ and evaluation map
$$\Ev:\ma\srarr (\proj^r)^n\times (\proj^r)^{r_0} \times (\proj^r)^{r_\infty}.$$
Let $l=\dim \vir{\ma}$. Then for $c\in
H^l((\proj^r)^{n+r_0+r_\infty})$, there exists a product of linear
subspaces $K$, Poincare-dual to $c$ so that $\Ev^{-1}(K)$ is a
finite union of reduced points, supported away from curves with
singular domain.}
\end{theorem}

Now, we specify the rubber moduli space we are considering. Let
genus be $0$.  Fix a degree $dL\in H_2(\proj^r)>0$, a number of
interior marked points $n$ and boundary marked points with
tangencies to $D_0$ and $D_\infty$, $r_0$ and $r_\infty$
respectively. Choose multiplicities $m_1^0\dots m_{r_0}^0$ and
$m_1^\infty\dots m_{r_\infty}^\infty$ so that
$(m_1^0+\dots+m_{r_0}^0)-(m_1^\infty+\dots+m_{r_\infty}^\infty)=md$.
Then, given a stable map $(f,C)\in
\cmbar_{0,n+r_0+r_\infty}(\proj^r,d)$ with smooth domain and
marked points,
$$\{x_1,\dots,x_n,p_1^0,\dots,p_{r_0}^0,p_1^\infty,\dots,p_{r_\infty}^\infty\}$$
the invertible sheaf
$$\cl=f^*L \otimes \oh(-(m_1^0p_1^0+\dots+m_{r_0}^0p_{r_0}^0)+(m_1^\infty
p_1^\infty+\dots+m_{r_\infty}^\infty p_{r_\infty}^\infty))$$ has
degree $0$.  $\cl$'s nonzero section, defined up to $\com^*$
induces a map
$$C \srarr P=\proj_{\proj^r}(L\oplus 1_{\proj^r})$$
giving a point in $\ma$.  Moreover, the automorphisms of the map
in $\cmbar_{0,n+r_0+r_\infty}(\proj^r,d)$ are in bijective
correspondence with the automorphisms of map in $\ma$.

\begin{theorem} {With the data on $\ma$ as above, $l=\dim\vir{\ma}$,
$c\in H^l((\proj^r)^{n+r_0+r_\infty})$, we have the following
equality of rubber and Gromov-Witten invariants
$$\deg(\Ev_{\cmbar}^*(c)\cap
\cmbar_{0,N}(\proj^r,d))=\deg(\Ev_\ma^*(c)\cap\vir{\ma}).$$}
\end{theorem}

\begin{proof} Pick a representative cycle $K$ as above.  Then
$\Ev_{\cmbar}^{-1}(K)\subseteq \cmbar_{0,N}(\proj^r,d)$ is a zero
dimensional reduced substack corresponding to maps with smooth
domains. $\Ev_{\ma}^{-1}(K)=\ft^{-1}\Ev_{\cmbar}^{-1}(K)$. By the
above consideration, given an integral zero dimensional substack,
$x$ in $\Ev_{\cmbar}^{-1}(K)$, $\ft^{-1}(x)$ is an integral zero
dimensional substack with the same automorphism group as $x$.

Consider the fiber square
$$\xymatrix{
D \ar[r]^>>>>{i'} \ar[d]& \cmbar_{0,n+r_0+r_\infty}(\proj^r,d)\ar[d]\\
K \ar[r]_<<<<<{i}  & (\proj^r)^{n+r_0+r_\infty} }$$ Now, the
refined Gysin map $i^!:A_j(\cmbar_{0,N}(\proj^r,d))\srarr
A_{j-l}(D)$ satisfies
$$\deg(i^!\vir{\cmbar_{0,N}(\proj^r,d)})=\deg(\Ev^*(c) \cap
\vir{\cmbar_{0,N}(\proj^r,d)}).$$ But since $K$ is a regularly
embedded substack, $i^!=(i')^*$.

Now, we need the following fact that ensures the compatibility of
the Gysin map with the virtual cycle construction (\cite{LT},
3.9).

\begin{lemma} {Let $\xi:X_0\srarr X$ be a regular embedding of
codimension $d$, $W$, a scheme, $W_0$, a scheme defined by the
following fiber square
$$\xymatrix{
W_0 \ar[r] \ar[d]& W\ar[d]\\
X_0 \ar[r]_\xi & X\\
}$$ If the tangent obstruction complexes $\ct_{W_0}^\bullet$ and
$\ct_W^\bullet$ are compatible in the sense of (\cite{LT}, 3.8)
and certain technical assumptions are satisfied, then
$$\xi^!\vir{W}=\vir{W_0}$$}
\end{lemma}

It can be checked that the induced virtual cycle on $D$ is just
$\vir{D}=[D]$.  Now recall (\cite{Vi}, 1.17) for a map from a
integral zero-dimensional stack $F$ to a point, $\pt$ (all stacks
over $\com$),
$$\deg(F/\pt)=\frac{1}{\delta(F)}$$
where $\delta(F)$ is the order of the automorphism group of $F$.

Therefore, the Gromov-Witten invariant is
$$\deg((i')^*\vir{\cmbar_{0,N}(\proj^r,d)})=\deg(D/\pt)=\sum_{F\in
D} \frac{1}{\delta(F)}.$$

Likewise, if $E$ is defined by
$$\begin{array}{ccc}
E & \srarr & \ma\\
\downarrow & & \downarrow \\
C & \srarr & (\proj^r)^N \\
& i &\\
\end{array}$$
the rubber invariant is
$$\deg((i')^*\vir{\ma})=\deg(E/\pt)=\sum_{F\in
E} \frac{1}{\delta(F)}.$$ \end{proof}

This result can be stated with beautiful succinctness following
(\cite{EGH} 2.9.2).  Consider rubber invariants into
$P=\proj_{\proj^r}(\oh(m)\oplus \oh)$. The numerical condition for
multiplicities implies
$$(m_1^0+\dots+m_{r_0}^0)-(m_1^\infty+\dots+m_{r_\infty}^\infty)=md.$$
Pick a homogeneous basis $\{a_1,\dots,a_v\}$ of $H^*(X)$.  Pick
variables $\beta_i$ dual to $a_i$ as in Definition \ref{accpot}.
Let $p_{k,i}$ be the variables corresponding to contact to $D_0$
with multiplicity $k$, dual to $a_i$.  $q_{k,i}$ corresponds
analogously to contact to $D_\infty$.  Consider a real variable
$x$. Let
$$P_j=\sum_{k\geq 1} p_{k,j}e^{-ikx}$$
$$Q_j=\sum_{k\geq 1} q_{k,j}e^{ikx}$$
Let $f$ be the rational Gromov-Witten potential of $\proj^r$, that
is
\begin{eqnarray*}
f(t_1,\dots,t_v,z)&=&\sum_d \sum_{n_i} \frac{1}{n_1!\dots n_v!}\\
&&\ \ ((\Ev_{\cmbar}^*((t_1a_1)^{n_1}\times\dots\times
(t_va_v)^{n_v}))\\
&&\ \cap \vir{\cmbar_{0,\sum n_i}(\proj^r,d)}).
\end{eqnarray*}

\begin{corollary} \label{fourierrat} {The rational rubber potential
is given by
$$A=\frac{1}{2\pi}\int_{0}^{2\pi} f(\beta_1+P_1+Q_1,\beta_2+P_2+Q_2,\dots,\beta_v+P_v+Q_v,ze^{imx})\ dx$$
where within the above formula, we treat $p_{k,i}$, $q_{k',i'}$ as
super-commuting variables of degree
$$\deg p_{k,i}=2-2\deg a_i-2k$$
$$\deg q_{k',i'}=2-2\deg a_{i'}+2k$$
and the $p$-variables are to be written before the $q$-variables.
}
\end{corollary}

\begin{proof} The operation
$$\frac{1}{2\pi}\int_{0}^{2\pi}\ dx$$
has the effect neglecting all coefficient of $e^{mix}$ for $m\neq
0$ which ensures that numerical condition for multiplicities is
satisfied. It is clear then that the rubber invariant is equal to
the corresponding Gromov-Witten invariant where interior and
boundary marked points are treated as marked points.
\end{proof}

This Fourier series formalism is similar to the residue formalism
of \cite{G}.

\subsection{Caporaso-Harris formula}

Here we examine rubber invariants of $\oh(1)$ over $\proj^1$
without powers of $\Top$ or $\Bot$.

\begin{lemma} All higher genus ($g\geq 1$) rubber invariants of
$L=\oh(1)$ over $X=\proj^1$ with at least one point of tangency to
$D_\infty$ and one interior marked point vanish.
\end{lemma}

\begin{proof}
%%Let us fix a degree $d$, a genus $g$, $m$ interior marked points, $r_0+r_\infty$ boundary marked points.  By the
%%dimension formula,
%%\begin{eqnarray*}
%%\vdim\ma&=&g-1+r_0+r_\infty+m+2d-1
%%\end{eqnarray*}
%%to get a non-vanishing invariants, we must have
%%$$g+2(d-1)\leq 0.$$

%%If $g=2$ then $d=0$ and we have a curve in a fiber of
%%$\proj_{\proj^1}(\oh(1)\oplus 1)\srarr \proj^1$.  Furthermore,
%%$\dim\ma=r_0+r_\infty+m$.  Therefore, to get a top dimensional
%%intersection, we must pull back a class from $H^2(\proj^1)$ at
%%every marked point, but this gives $0$ for the intersection. If
%%$g=1$ then $d=0$ and $\vdim\ma=r_0+r_\infty+m-1$.  The curve must
%%lie in a fiber, and so the invariant, again must vanish.
This is a virtual dimension count.
\end{proof}

Now, following \cite{EGH}, let us apply Corollary \ref{fourierrat}
to compute the rational rubber potential.  Let us change notation
slightly and write a basis for $H^*(\proj^1)$ as
$$a_0\in H^0(\proj^1)$$
$$a_2\in H^2(\proj^1)$$
and write $t_0,t_2$ for variables dual to $a_0, a_2$.  The
rational Gromov-Witten potential for $\proj^1$ is
$$f_{\proj^1}(t_0,t_2,z)=\frac{t_0^2 t_2}{2}+e^{t_2} z.$$
Therefore, the rubber potential, $A$ satisfies

%%is
%%\begin{eqnarray*}
%%\hbar^{-1}\frac{1}{2\pi}\int_0^{2\pi} \frac{(\beta_0+\sum_k
%%p_{k,0}e^{-ikx}+\sum_k q_{k,0}e^{ikx})^2(\beta_2+ \sum_k
%%p_{k,2}e^{-ikx}+\sum_k q_{k,2}e^{ikx})}{2} dx\\
%%+\hbar^{-1}\frac{1}{2\pi}\int_0^{2\pi}
%%e^{\beta_2+\sum_k p_{k,2}e^{-ikx}+\sum_k q_{k,2}e^{ikx}} ze^{ix}
%%dx
%%\end{eqnarray*}

\begin{eqnarray*}
A_2&\equiv &\frac{\partial A}{\partial
\beta_2}|_{\beta_2=0}\\
&=&\hbar^{-1}\left(\beta_0^2+\sum_{k}
p_{k,0}q_{k,0}+\frac{1}{2\pi}\int_0^{2\pi} e^{\sum_k
p_{k,2}e^{-ikx}+\sum_k q_{k,2}e^{ikx}+ix} z\ dx\right).
\end{eqnarray*}

Let us write down the relative potential of $(\proj^2,L)$, that is
the projective plane relative a line.  Let us choose
$\{H^2,H^1,1\}$ as a basis of $H^*(\proj^2)$.  We restrict
ourselves to the potential involving only cohomology of the form
$H^2$ and at least one $p$ variable.  Let us use $\theta_1$ to
express the element of $\cf$ dual to $H^2$. Let us use $p_{n,0}$
and $p_{n,2}$ to express $n$th order multiplicities to $H$ at $1$
and $[\text{pt}]\in H^*(L)$ respectively. Let us write the degree
as $z^d$ where $d$ denotes the class of $dL\in H_2(\proj^2)$.
Therefore, $F$ is an expression in $\thbar$, $\theta_1$,
$\tp_{n,0}$, $\tp_{n,2}$, and $z^d$.

%%Note that $F|_{\theta_1=0}$ counts relative maps in $(\proj^2,L)$
%%ith no interior marked points.  If we have $r$ boundary marked
%%points,
%%$$\vdim\mz=(2-3)(1-g)+<f^*c_1(T\proj^2)-f^*L,dL>+r=(g-1)+2d+r$$
%%Since the highest dimensional class we can pull back at a boundary
%%marked point is of dimension $2$, the only way to get a top
%%dimensional intersection is if $d=0$ and $g=1$ which we have
%%decided to exclude.  Therefore $F|_{\theta_1=0}=0$.

By dimensional considerations, $F|_{\theta_1=0}=0$ and the
differential equation of Theorem \ref{Freconstruct} becomes
$$\frac{\partial F}{\partial \theta_1}=A_2\cdot F.$$
Unwinding the action of $\ccr$ on $\cf$, we see that this becomes
$$\frac{\partial F}{\partial \theta_1}=
(\sum_{k} k\tp_{k,0}\frac{\partial}{\partial
\tp_{k,2}}+\hbar^{-1}\frac{1}{2\pi}\int_0^{2\pi} e^{\sum_k
\tp_{k,2}e^{-ikx}+\sum_k k\thbar\frac{\partial}{\partial
\tp_{k,0}}e^{ikx}+ix} z\ dx)\ F.$$

This is the expression of the Caporaso-Harris formula as written
in \cite{G}.

\subsection{Ruled Surfaces $\eff_n$}

We can apply the rubber formalism to derive the inductive formula
for the relative Gromov-Witten invariants on Hirzebruch surfaces
from \cite{Va2}. Let $\eff_n$ be the ruled surface
$$\eff_n=\proj_{\proj^1}(\oh(n)\oplus \oh)$$
where $n\geq 0$. Let $\pi$ be the projection
$\pi:\eff_n\srarr\proj^1$

Let $D\subset \eff_n$ be the infinity section of
$\oh_{\proj^1}(n)$  The second homology of $\eff_m$, $H_2(\eff_m)$
is generated by $C_0=D$, and $f$, a fiber of $\pi$.

%%$C_0\cdot C_0=-m$, $C_0\cdot f=1$, $f\cdot f=0$.

%%Note that we have
%%$$c_1(T\eff_m)=(2+l)f+2C_0$$

%%Let us imitate the above section and look at the potential where
%%the only cohomology we allow at interior marked points are point
%%classes. The virtual dimension of the stack $\mz$ of maps of
%%degree $aC_0+bf\in H_2(\eff_m)$ is
%%$$\vdim \mz=g-1+2a+b+r+n$$

%%Now the effective cone is generated by $C_0$ and $f$.  Therefore,
%%we must have $a\geq 0$, $b\geq 0$. To get a zero dimensional
%%intersection on the moduli stack
%%$$g-1+2a+b\leq 0$$

Again, let us consider the terms in $F$ with at least one
$p$-variable and no point classes at interior marked points.  By
dimensional reasons, the only non-vanishing invariant comes from
degree 1 maps to a fiber. In this case, the virtual cycle of the
moduli space agrees with the usual fundamental cycle, yielding

%%If $g=1$, then we have only constant curves which we've decided to
%%exclude. If $g=0$, we must have $b=1$, so our curve is a fiber.
%%Therefore, our moduli space is a copy of $\proj^1$ where we have
%%exactly one boundary marked point. It is easy to see by writing
%%down the tangent/obstruction complex that the virtual cycle agrees
%%with the usual fundamental class. Therefore, we get a nonzero
%%invariant from pulling back a point class on $D=\proj^1$ at the
%%boundary marked point. Therefore,
$$F|_{\theta_1=0}=\tp_{1,2} z^f.$$

Now, we can compute the rubber invariants of $X=\proj^1$,
$L=N_{D/\eff_m}=\oh(-m)$.  By the same arguments as above,
$$A_2=\hbar^{-1}\left(\beta_0^2+\sum_{k} p_{k,0}q_{k,0}+\frac{1}{2\pi}\int_0^{2\pi}
e^{\sum_k p_{k,2}e^{-ikx}+\sum_k q_{k,2}e^{ikx}-imx} z\
dx\right).$$

We get the differential equation
$$\frac{\partial F}{\partial \theta_1}=
(\sum_{k} k\tp_{k,0}\frac{\partial}{\partial
\tp_{k,2}}+\hbar^{-1}\frac{1}{2\pi}\int_0^{2\pi} e^{\sum_k
\tp_{k,2}e^{-ikx}+\sum_k k\thbar\frac{\partial}{\partial
\tp_{k,0}}e^{ikx}-imx} z^{C_0}\ dx)\ F.$$

Under the identification
$$\Bl_{0} \proj^2=\eff_1$$
the above recursion formula reduces to Ran's \cite{Ra1}.

\subsection{Rational Gromov-Witten Invariants for $\proj^n$}

Here we consider the rational relative Gromov-Witten invariants of
$(\proj^n,H)$ where $H$ is a hyperplane in $\proj^n$ and $n\geq
3$.  We follow the computation of \cite{EGH} which gives the
formula of \cite{Va1}. Let us compute the potential $F$ where we
count positive degree curves and at interior marked points we pull
back cohomology classes from $H^k(\proj^n)$ where $k\geq 4$.
Consider cohomology classes
$$e_4,e_6,\dots,e_{2n}$$
where $e_{2i}=[H]^{i}$.  Let $\theta_{2i}$ be dual to $e_{2i}$ On
$H^*(H)=H^*(\proj^{n-1})$, we look at cohomology classes
$$c_0,c_2,\dots,c_{2n-2}$$
where $c_{2i}$ generates $H^{2i}(\proj^r)$.
Let $\tp_{k,i}$ be dual to $c_{2i}$.

%%If we have $r$ boundary marked points and $m$ interior marked
%%points and are looking at curves of degree $d$, then
%%$$\vdim\mz=n(d+1)+r+m-3$$

%%We will induct on the number of interior marked points.  Let us
%%find the seed values for $m=0$. We can pull back classes of
%%codimension up to $2n-2$ at boundary marked points. Therefore, to
%%get a top dimensional intersection, we need
%%$$n(d+1)+r-3-r(n-1)\leq 0$$
%%Since $r\leq C\cdot H=d$, we must have
%%$$2d+n-3\leq 0$$
%%which is impossible unless $d=0$ and we chose as our convention
%%not to count constant curves.

By dimensional considerations,
$$F|_{\theta_4=\dots=\theta_{2n}=0}=0.$$

Let
$$f_{\proj^{n-1}}(t_0,t_2,\dots,t_{2n-2},z)$$
be the rational Gromov-Witten potential of $\proj^{n-1}$ where
$t_{2i}$ is dual to a cohomology class in $H^{2i}(\proj^{n-1})$
and $z$ is the degree marker.  Then, we can use Corollary
\ref{fourierrat} to write the rubber potential,
$$A_{2i}=\frac{1}{2\pi}\int_{0}^{2\pi}
\frac{\partial f}{\partial
t_{2i}}(\beta_1+P_1+Q_1,\beta_2+P_2+Q_2,\dots,\beta_v+P_v+Q_v,ze^{ix})\
dx.$$

Our differential equations become
$$\frac{\partial F}{\partial \theta_{2i}}=A_{2i-2} \cdot F.$$

%\bibliographystyle{plain}
%\bibliography{thesis}

%\end{document}

%% file: hamiltonian.tex
In this section, we study a formalism for relative Gromov-Witten
invariants that duplicates the structure of the Symplectic Field
Theory \cite{EGH} of Eliashberg, Givental, and Hofer.  This
formalism has the added advantage that it takes into account some
of the redundancies of rubber invariants given by Theorem
\ref{tcbsc}.  The rubber invariants are encoded in a certain
generating function called the {\em Hamiltonian}

Let $X$ be a projective manifold and let $P=\proj(L \oplus 1_X)$
be the projective completion in a line-bundle $L$ over $X$.  We
considered $\ma$, a moduli space of curves in $P$ relative the
zero and infinity sections and quotiented by a $\cs$-action that
dilates the fibers.  This moduli space possesses an evaluation map
$$\Ev:\ma \srarr X^m \times X^{r_0} \times X^{r_\infty}$$
where the three factors in the target denote the image of $n$
interior marked points, $r_0$ boundary marked points evaluating to
the zero section and $r_\infty$ boundary marked points evaluating
to the infinity section.

In Symplectic Field Theory, there is a similar moduli space,
$\cm$.  The construction and compactification of this moduli space
are markedly different and the evaluation map is
$$\Ev:\cm \srarr (S(L))^m \times X^{r_0} \times X^{r_\infty}$$
where $S(L)$ is the unit circle bundle in $L$.  Consequently, the
classes that are pulled back at interior marked points are from
$H^*(S(L))$. Now the cohomology of $S(L)$ is related to that of
$X$ by the following Gysin sequence.

$$\xymatrix{
H^*(L,L_0) \ar[r] \ar@{=}[d] & H^*(L) \ar[r]\ar@{=}[d] &
H^*(L_0)\ar[r]\ar@{=}[d] & H^{*+1}(L,L_0) \ar[r]\ar@{=}[d]& H^{*+1}(L)\ar@{=}[d]\\
H^{*-2}(X) \ar[r]_>>>>{\cup c_1(L)} & H^*(X) \ar[r]_>>>>{\pi^*}&
H^*(S^1(L)) \ar[r]_>>>>>{\pi_*} & H^{*-1}(X)\ar[r]_{\cup c_1(L)}&
H^{*+1}(X)}$$ where $\pi:S(L) \srarr X$.

Classes in $H^*(S^1(L))$ are non-canonically decomposable into
classes of two kinds, those in the image of $\pi^*$ and those in
the cokernel of $\pi^*$. Classes in the image of $\pi^*$ do not
fix the $S^1$ symmetry of the fibers in $S^1(L)$ while classes in
the cokernel of $\pi^*$ do.  Classes in the cokernel, we consider
{\em phase-fixing} marked points while those in the image we
consider {\em non-phase-fixing}. There is a subset of Symplectic
Field Theory invariants that involve pulling back phase-fixing
classes at exactly one marked point.  It is this theory that we
translate into our context.  As an aside, we note that we believe
that the analog of SFT in which we use a number of phase-fixing
classes is the invariants as below, but enriched by powers of
$c_1(\Bot)$.

Consider two types of marked point depending on the type of
cohomology classes that is pulled back along their evaluation
maps.

\begin{enumerate}
\item[1)] Phase-fixing marked points $\phi$ that involve a class
of the form $\ev_\phi^* (\pi_*a)$ or alternatively by a class
$\ev_\phi^* b$ where $b \in \ker(\cup c_1(L) : H^*(X) \srarr
H^{*+2}(X))$.

\item[2)] Non-phase fixing marked points $i$ that involve $ev_i^*(a)$ where
$a$ is arbitrary but can be considered to be a representative of a class
in $\coker(\cup c_1(L) : H^{*-2}(X) \srarr H^*(X))$.
\end{enumerate}

\subsection{Definition of Hamiltonian}

The Hamiltonian will be dependent on the choice of a particular
kind of basis.

\begin{definition} {A {\em phase fixing basis} is a homogeneous basis
$\{a_1,\dots,a_m\}$ for $\ker(\cup c_1(L):H^*(X)\srarr
H^{*+2}(X))$}
\end{definition}

\begin{definition} {A {\em non-phase fixing basis} is a set of
homogeneous elements $\{b_1,\dots,b_n\}$ of $H^*(X)$ that descend
to a basis of $\coker(\cup c_1(L):H^{*-2}(X)\srarr H^*(X))$.}
\end{definition}

We define $\ch$, a graded super-Weyl algebra over $\rat[H_2(X)]$.
Pick a phase-fixing basis $\{a_1,\dots,a_m\}$, a non-phase fixing
basis $\{b_1,\dots,b_n\}$, and a homogeneous basis
$\{c_1,\dots,c_s\}$ of $H^*(X)$.  We will consider variables
$\sig{i},\tau_j$ corresponding to $a_i,b_j$ respectively.  For
every positive integer $n$, we will have $p_{n,i}, q_{n,i}$
corresponding to $c_i$.  Let us consider a graded algebra
generated by the following elements
\begin{eqnarray*}
\deg \hbar&=&-2(\dim X-2)\\
\deg z^A&=&2<c_1(TX),A>+<c_1(L,A)>\\
\deg \sig{i}&=&-1-\deg a_i\\
\deg \tau_i&=&2-\deg b_i\\
\deg p_{n,i}&=&2-\deg c_i-2n\\
\deg q_{n,i}&=&2-\deg c_i+2n.
\end{eqnarray*}

Note that we have changed the grading of $\sig{i}$ from what we
would have expected in the definition of the rubber potential. In
$\ch$, $\hbar$, $z^A$, $\sig{i}$, $\tau_i$ are supercentral while
the $p$- and $q$-variables satisfy
$$[p_{n_1,i_1},p_{n_2,i_2}]=0,\ [q_{n_1,i_1},q_{n_2,i_2}]=0,\
[q_{n_1,i_1},p_{n_2,i_2}]=n_1\delta_{n_1,n_2}g^{i_1 i_2}\hbar.$$

As in $\ccr$, the algebra $\ch$ consists of Laurent series in
$\hbar$ whose coefficients are power series in the $p$-variables
whose coefficients are power series in the $\tau$- variables whose
coefficients are polynomials in the $\sigma$ and $q$-variables.

To define the Hamiltonian, we define the following formal sums
\begin{eqnarray*}
S&=&\sum \sig{i} a_i,\\
T&=&\sum_n \frac{1}{n!}(\sum_i \tau_i b_i)^n,\\
P&=&\sum_n \frac{1}{n!}\left(\sum_{k,i} p_{k,i}\right)^n,\\
Q&=&\sum_n \frac{1}{n!}\left(\sum_{k,i} q_{k,i}\right).^n
\end{eqnarray*}

\begin{definition} \label{hampot} The {\em Hamiltonian} $H$ is
$$H=\sum_\chi \sum_{A\in H_2(X)}  \hbar^{-\frac{1}{2}\chi} (S\cdot T,P,Q)_{\chi,A,0} z^A$$
where $S\cdot T$ denotes multiplication in the algebra $\ch$.
\end{definition}

$H$ is linear in the $\sigma$-variables.  The $\sigma$ variable
keeps track of the phase-fixing classes which are evaluated at a
single, distinguished marked point. $H$ is homogeneous of degree
$-1$.

\subsection{Dependence on Representatives}

We will show how the Hamiltonian depends on the choice of
representative classes in $\coker(\cup c_1(L) : H^{*-2}(X) \srarr
H^*(X))$.

\begin{lemma} \label{Ldeped} Consider $\ma$ with two
distinguished marked points, $\phi$ and $i$. Let $a\in H^*(X)$
satisfy $a\cup c_1(L)=0$.  Then

\noindent
\setlength{\unitlength}{.5cm}
\begin{center}
\begin{picture}(11,5)
  \put(1,3)        {$\ev_\phi^*(a)\ev_i^*(c_1(L))=\ev_\phi^*(a)($}
  \put(9.5,3.15)    {\line(1,-3){.5}}
  \put(9.5,2.85)    {\line(1, 3){.5}}
  \put(9.6,2.2)   {\line(1,0) {.5}}
  \put(9.6,3.8)   {\line(1,0) {.5}}
  \put(10.2,2)     {$\phi$}
  \put(10.2,3.6)   {$i$}

  \put(10.5,3)       {$-$}

  \put(11.5,3.15)    {\line(1,-3){.5}}
  \put(11.5,2.85)    {\line(1, 3){.5}}
  \put(11.6,2.2)   {\line(1,0) {.5}}
  \put(11.6,3.8)   {\line(1,0) {.5}}
  \put(11.2,2)     {$i$}
  \put(11.2,3.6)   {$\phi$}

  \put(12,3)      {$)$}

\end{picture}
\end{center}
\end{lemma}

\begin{proof} Multiply the formula from Theorem \ref{tcbsc} with
$\ev_\phi^*(c_1(L))$. \end{proof}

If we view $\phi$ as the phase-fixing marked point, every term in
the Hamiltonian will involve a factor of the form $\ev_\phi^* a$
where $a\cup c_1(L)=0$

Let us pick a phase-fixing basis, $$\{a_1,\dots,a_m\}\subset
H^*(X)$$ and a non-phase fixing basis
$$\{b_1,\dots,b_n\}\subset H^*(X).$$  From
the construction of our algebra $\ch$, the Hamiltonian is
invariant under change of basis of the form
$$b_i'=M_{ij}b_j.$$
Therefore, we need only determine how $H$ varies when we change
$b_i\in H^*(X)$ representing a class $[b_i]\in\coker(\cup c_1)$.
Let us change $\{b_1,\dots,b_n\}$ one element at a time. Write
$$b_1'-b_1=dc_1(L).$$
Define $\tilde{b_t}=(1-t)b_1+tb_1'$.  Let $H_t$ be $H$ with
$\tilde{b_t}$ substituted for $b_1$ so $H_0$ is the Hamiltonian
with $b_1$ in place, $H_1$ is the Hamiltonian with $b_1'$ in
place.

$$H_t= \sum_\Gamma \frac{1}{
\Aut_\Gamma(R_0)\Aut_\Gamma(R_\infty) k_1! \dots k_n! }
<a_j,\tilde{b_t}^{k_1}b_2^{k_2}\dots b_n^{k_n}>$$
$$\cdot\hbar^{\frac{\chi}{2}}\sig{i}\tau_1^{k_1}\dots\tau_n^{k_n}p^\Gamma
q^\Delta z^d$$ so
$$\frac{\partial H_t}{\partial t}=
\sum_\Gamma
\frac{1}{\Aut_\Gamma(R_0)\Aut_\Gamma(R_\infty)k_1!\dots k_n!}$$
$$\cdot<a_j,k_1(dc_1(L))\tilde{b_t}^{k_1-1}b_2^{k_2}\dots b_n^{k_n}>
\hbar^{-\frac{1}{2}\chi}\sig{i}\tau_1^{k_1}\dots\tau_n^{k_n}p^\Gamma
q^\Delta z^d.$$

But
$$<a_j,(dc_1(L))\tilde{b_t}^{k_1-1}b_2^{k_2}\dots b_n^{k_n}>=$$
\noindent \setlength{\unitlength}{.5cm}
\begin{center}
\begin{picture}(11,5)
  \put(1,3)        {$\Ev^*(a_j\times d \times \tilde{b_t}^{k_1-1} b_2^{k_2}\dots
b_n^{k_n})\cdot($}
  \put(11.5,3.15)    {\line(1,-3){.5}}
  \put(11.5,2.85)    {\line(1, 3){.5}}
  \put(11.6,2.2)   {\line(1,0) {.5}}
  \put(11.6,3.8)   {\line(1,0) {.5}}
  \put(12.2,2)     {$\phi$}
  \put(12.2,3.6)   {$1$}

  \put(12.5,3)       {$-$}

  \put(13.5,3.15)    {\line(1,-3){.5}}
  \put(13.5,2.85)    {\line(1, 3){.5}}
  \put(13.6,2.2)   {\line(1,0) {.5}}
  \put(13.6,3.8)   {\line(1,0) {.5}}
  \put(13.2,2)     {$1$}
  \put(13.2,3.6)   {$\phi$}

  \put(14.5,3)      {$).$}

\end{picture}
\end{center}
where the first two entries in $\Ev$ corresponds to the marked
points denoted by $\phi$ and $1$ respectively.

Let us define the generating function $K_t$ by
$$K_t=\sum_\Gamma
\frac{1}{\Aut_{\Gamma}(R_0)\Aut_{\Gamma}(R_\infty)k_1!\dots
k_n!}$$
$$\cdot <d_1\tilde{b_t}^{k_1}b_2^{k_2}\dots b_n^{k_n}>
\hbar^{\frac{\chi}{2}}\tau_1^{k_1+1}\dots\tau_n^{k_n}p^\Gamma
q^\Delta z^d.$$

Note that $K$ is of even degree and so
$$\frac{\partial
H_t}{\partial t}=[K_t,H_t].$$

\subsection{The Hamiltonian as a Differential}

We may define a differential, $D^H$, on the algebra $\ch$.  One
can compute the homology with respect to this differential. This
homology will be an invariant of $(X,L)$.

\begin{lemma} \label{pfshift} Consider $\ma$ with at least two
marked points including $\phi_1$ and $\phi_2$.  Let
$a_1,a_2\in\ker(\cup c_1(L):H^*(X)\srarr H^{*+2}(X))$. Then

\noindent \setlength{\unitlength}{.5cm}
\begin{center}
\begin{picture}(16,5)
  \put(1,3)        {$\ev_{\phi_1}^*(a_1)\ev_{\phi_2}^*(a_2)$}
  \put(7,3.15)    {\line(1,-3){.5}}
  \put(7,2.85)    {\line(1, 3){.5}}
  \put(7.1,2.2)   {\line(1,0) {.5}}
  \put(7.1,3.8)   {\line(1,0) {.5}}
  \put(6.4,2)     {$\phi_1$}
  \put(6.4,3.6)   {$\phi_2$}

  \put(8,3)       {$=\ev_{\phi_1}^*(a_1)\ev_{\phi_2}^*(a_2)$}

  \put(14.9,3.15)    {\line(1,-3){.5}}
  \put(14.9,2.85)    {\line(1, 3){.5}}
  \put(15,2.2)   {\line(1,0) {.5}}
  \put(15,3.8)   {\line(1,0) {.5}}
  \put(14.3,2)     {$\phi_2$}
  \put(14.3,3.6)   {$\phi_1$}

\end{picture}
\end{center}
\end{lemma}

\begin{proof} We intersect

\noindent
\begin{center}
\setlength{\unitlength}{.5cm}
\begin{picture}(17,5)
  \put(1,3)  {$\ev_{\phi_1}^*(c_1(L))-\ev_{\phi_2}^*(c_1(L))$}
  \put(9.7,3)  {$=$}
  \put(12,3.15){\line(1,-3){.5}}
  \put(12,2.85){\line(1, 3){.5}}
  \put(12.1,2.2) {\line(1,0) {.5}}
  \put(12.1,3.8) {\line(1,0) {.5}}
  \put(11.4,2)   {$\phi_1$}
  \put(11.4,3.6)   {$\phi_2$}
  \put(13.7,3)   {$-$}
  \put(16,3.15){\line(1,-3){.5}}
  \put(16,2.85){\line(1, 3){.5}}
  \put(16.1,2.2) {\line(1,0) {.5}}
  \put(16.1,3.8) {\line(1,0) {.5}}
  \put(15.4,2)   {$\phi_2$}
  \put(15.4,3.6)   {$\phi_1$}
\end{picture}
\end{center}

\noindent with $\ev_{\phi_1}^*(a_1)\ev_{\phi_2}^*(a_2)$.
\end{proof}

\begin{theorem} {$H$ is nilpotent of order $2$ in $\ch$.}
\end{theorem}

\begin{proof} Write
$$H=\sum_i \sig{i} H^i.$$
Notice that $H^i$ is an element of degree $1-\deg\sig{i}=2+\deg
a_i$.
\begin{eqnarray*}
H^2&=&(\sum_i \sig{i} H^i)^2\\
&=&\sum_i
\sig{i}H^i\sig{i}H^i+\sum_{i<j}(\sig{i}H^i\sig{j}H^j+\sig{j}H^j\sig{i}H^i)\\
&=&\sum_i \sig{i}^2 (H^i)^2+\sum_{i<j} (\sig{i}\sig{j}((-1)^{(\deg
a_i)(\deg a_j-1)}H^i H^j\\
&&\ \ +(-1)^{(\deg a_i-1)(\deg a_j)+(\deg a_i-1)(\deg
a_j-1)}H_jH_i)).
\end{eqnarray*}

Note that if $\deg a_i$ is odd, then $\sig{i}$ is even, so
$\sig{i}^2\neq 0$, so we must show $(H^i)^2=0$.  But by arguments
analogous to the proof of the degeneration formulae, $(H^i)^2$ is
given by \noindent
\begin{center}
\setlength{\unitlength}{.5cm}
\begin{picture}(17,5)
  \put(1,3)        {$H_i^2=(\ev_{\phi_1}^*(a_1)\ev_{\phi_2}^*(a_2)(\text{sum over terms})\cup$}
  \put(15.2,3.15)  {\line(1,-3){.5}}
  \put(15.2,2.85)  {\line(1, 3){.5}}
  \put(15.3,2.2)   {\line(1,0) {.5}}
  \put(15.3,3.8)   {\line(1,0) {.5}}
  \put(14.6,2)     {$\phi_1$}
  \put(14.6,3.6)   {$\phi_2$}
  \put(15.9,3)       {$)\cap\ \vir{\ma}.$}
\end{picture}
\end{center}
where "sum over terms" indicate a sum over all the cohomology
classes and formal variables that enter into Definition
\ref{hampot}. But by the above lemma, this equals
\begin{center}
\setlength{\unitlength}{.5cm}
\begin{picture}(24,5)
  \put(1,3)        {$(\ev_{\phi_1}^*(a_i)\ev_{\phi_2}^*(a_i)(\text{sum over terms})\cup$}
  \put(12.9,3.15)  {\line(1,-3){.5}}
  \put(12.9,2.85)  {\line(1, 3){.5}}
  \put(13,2.2)   {\line(1,0) {.5}}
  \put(13,3.8)   {\line(1,0) {.5}}
  \put(12.3,2)     {$\phi_2$}
  \put(12.3,3.6)   {$\phi_1$}
  \put(13.6,3)     {$)\cap\ \vir{\ma}$}
\end{picture}
\begin{picture}(24,3)
  \put(1,3)        {$=-(\ev_{\phi_2}^*(a_i)\ev_{\phi_1}^*(a_i)(\text{sum over terms})\cup$}
  \put(14.3,3.15)  {\line(1,-3){.5}}
  \put(14.3,2.85)  {\line(1, 3){.5}}
  \put(14.4,2.2)   {\line(1,0) {.5}}
  \put(14.4,3.8)   {\line(1,0) {.5}}
  \put(13.7,2)     {$\phi_1$}
  \put(13.7,3.6)   {$\phi_2$}
  \put(15.1,3)     {$)\cap\ \vir{\ma}=-H_i^2.$}
\end{picture}
\end{center}

Since $\sig{i}\sig{j}\neq 0$, we must show
\begin{eqnarray*}
0&=&(-1)^{(\deg a_i)(\deg a_j-1)}H^i H^j+(-1)^{(\deg a_i-1)(\deg
a_j)+(\deg a_i-1)(\deg a_j-1)}H^j H^i\\
&=&(-1)^{\deg a_i-1}(H^j H^i-(-1)^{(\deg a_i)(\deg a_j)}H^i H^j).
\end{eqnarray*}
The proof is analogous to the one above. \end{proof}

\begin{definition} {The operator $D^H$, $D^H:\ch\srarr\ch$ is
defined on homogeneous elements of $\ch$ by $D^H:f\mapsto
Hf-(-1)^{\deg f} fH$.}
\end{definition}

\begin{corollary} {$D^H$ is a differential on the graded algebra
$\ch$}\end{corollary}

\begin{proof} It is a simple computation to show $(D^H)^2=0$ and
$D^H(fg)=D^H(f)g+(-1)^{\deg f}f D^H(g)$.\end{proof}

\subsection{Hamiltonian Homology}

The chain complex $(\ch,D^H)$ is invariant under changing the
choice of representatives of $\coker(\cup c_1(L))$ in the
non-phase fixing basis.

\begin{theorem} {Let $\tilde{b_t}=(1-t)b_1+tb_1'$.  Let $H_t$ be
$H$ with $\tilde{b_t}$ substituted for $b_1$.  Then there is a
family of isomorphisms $Q_t:\ch\srarr \ch$ so that
$$D^{H_t} Q_t(f)=Q_t(D^{H}f).$$}
\end{theorem}

\begin{proof}

Let $K_t$ be as above.  Define $A_t$ by
$$A_t=\int_0^t K_s\ ds,$$
and let $Q_t$ be given by
$$Q_t(f)=e^{A_t} f e^{-A_t}.$$
Write $f_t=Q_t(f)$ and let
$$E_t=D^{H_t}(f_t)-Q_t(D^{H_0}f).$$
Therefore,
$$E_t=[H_t,f_t]-Q_t([H_0,f],$$
\begin{eqnarray*}
\frac{\partial E_t}{\partial
t}&=&[[K_t,H_t],f_t]+[H_t,[K_t,f_t]]-[K_t,Q_t([H,f])]\\
&=&[K_t,E_t].
\end{eqnarray*}
Since $E_0=0$, it follows that $E_t=0$. \end{proof}

By changing representatives of $\coker(\cup c_1(L))$ one-by-one,
we see that the differential graded algebra and hence the
Hamiltonian homology is invariant.

There are other versions of the Hamiltonian Homology, rational
Hamiltonian homology and contact Hamiltonian homology.  They bear
the same relation to Hamiltonian Homology as their analogs do to
Symplectic Field Theory Homology.  We refer the reader to
\cite{EGH} for details.
%\noindent
%\bibliographystyle{plain}
%\bibliography{thesis}

%\end{document}

%% file: local.tex
In this section, we give a proof of the degeneration formula in
Theorem \ref{degenma}.  That degeneration formula encodes in a
generating function relation (3) among line-bundles on $\ma$:
$$\Bot\otimes\ev_i^*L=\Lnb{i}.$$

A proof is outlined above, but here we give a more direct proof
using the virtual localization technique from \cite{K} and
\cite{GP} and adapted for the relative case in \cite{GV}.

Our strategy is to evaluate the equivariant cap product,
$$\hbar \cup \Ev^*c \cup \ev_i^*(c_1(\oh(1))+c_1(L))) \cap \vir{\my}.$$
which we know to be zero on a particular stack, $\my$ by
localization. The localization formula will give a relation among
cycle classes which when intersected with cohomology classes will
give Theorem \ref{degenma}.

\subsection{Target Schemes}

We need to construct a stack $\my$ that is closely related to
$\ma$.  Let $X$ be a projective manifold and $L$ be a line-bundle
over $X$.  Let $P=\proj_X(L\oplus 1_X)$. Let $p:P\srarr X$ be the
projection.  Let $i_0:D_0\srarr P$, $i_\infty:D_\infty\srarr P$ be
the inclusions of the zero and infinity sections respectively.

We want to consider stable maps into $P$ relative $D_0$ and
$D_\infty$.  The stack of stable maps we construct will differ
from $\ma$ in that we do not have a $\com^*$-action that dilates
the fiber of $P$.  The construction, however, is analogous to that
of $\mz$ and $\ma$.

Define the scheme $\ _{k,l} Y$ as the union of $k+l+1$ copies of
$P$,
$$\ _{k,l} Y=P_{-k}\sqcup_X \dots \sqcup_X P_{-1} \sqcup_X P_0 \sqcup_X P_1
\sqcup_X \dots \sqcup_X P_l$$ where $X_0\subset P_i$ is identified
with $X_\infty\subset P_{i+1}$.  Let the automorphism group of $\
_{k,l} Y$ be $(\cs)^k\times (\cs)^l$ where the first $k$ copies of
$\cs$ dilate the fibers of $P_{-k},\dots,P_{-1}$ and the last $l$
copies of $\cs$ dilate the fibers of $P_1,\dots,P_l$.  Note that
there is no $\cs$-factor dilating $P_0$.  $\my$ is the stack of
stable pre-deformable maps to $\ _{k,l} Y$ with data given by a
rubber graph.

The rigorous definition of $\my$ is analogous to those of $\ma$
and $\mz$. Begin by defining triples
$(Y[k,l],D_0[k,l],D_\infty[k,l])$ indexed by a pair of
non-negative integers where $Y[k,l]$ is a projective manifold with
a $G[m]\times G[n]$-action, and $D_0[k,l]$ and $D_\infty[k,l]$ are
smooth divisors. Let
\begin{eqnarray*}
Y[0,0]&=&P \\
D_0[0,0]&=&D_0\\
D_\infty[0,0]&=&D_\infty
\end{eqnarray*}
where $D_0$ and $D_\infty$ are the zero and infinity sections in
$P$. We define $Y[k,l]$ inductively,
$$Y[k+1,l]=\Bl_{D_\infty[k,l]\times\{0\}} (Y[k,l]\times\aff^1).$$
$D_\infty[k+1,l]$ is the proper transform of
$D_\infty[k,l]\times\aff^1$, $D_0[k+1,l]$ is the inverse image of
$D_0[k,l]\times\aff^1$.
$$Y[k,l+1]=\Bl_{D_0[k,l]\times\{0\}} (Y[k,l]\times\aff^1).$$
$D_0[k+1,l]$ is the proper transform of $D_0[k,l]\times\aff^1$,
$D_\infty[k+1,l]$ is the inverse image of
$D_\infty[k+1,l]\times\aff^1$. The $G[k]\times G[l]$ actions
dilates the fibers in the tails.

Given a rubber graph $\Gamma$, we construct $\my=\cm(\cy,\Gamma)$
by mimicking the construction of $\mz$ and $\ma$.  We consider
families of pre-deformable relative maps described by $\Gamma$
that are stable under the $G[m]\times G[n]$-action.  We glue these
families together into a stack and then quotient by the
$G[m]\times G[n]$-action.   $\my$ carries a virtual cycle.

Note that $\my$ is different from $\mz(P,D_0\sqcup D_\infty)$
since in the construction of $Z[n]$, $D_0$ and $D_\infty$ are
blown up simultaneously and the $\com^*$-action dilates the fibers
of their exceptional divisors simultaneously.

\subsection{Equivariant Data}

We will perform a virtual localization computation on $\my$.  We
first define a $\com^*$ action on $P$.  Specify points of
$P=\proj_X(L\oplus 1)$ by $[l:t]$.  For $\lambda\in\com^*$, define
the group action by
$$\lambda\cdot [l:t]=[\lambda l:t].$$
This $\com^*$ action induces a $\com^*$-action on $\my$ so that
$$\Ev:\my\srarr P^m\times X^{r_0}\times X^{r_\infty}$$
is equivariant.

Now, write $\HS{*}(\pt)=\com[\hbar]$ where $\hbar$ is the Euler
class of the equivariant line bundle on $\pt$ under the group
action
$$\lambda \cdot t = \lambda t.$$

Now, let $\Gamma$ be a rubber graph where there is at least one
vertex that does not correspond to a trivial cylinder. The proof
of the virtual localization theorem holds for $\cm(\cy,\Gamma)$
with trivial modifications.

Let $\pi:P\srarr X$.  Consider the composition
$$\Ev_X=(\pi^m\times \id_{X^{r_0}} \times \id_{X^r_\infty}):\my\srarr P^m\times
X^{r_0}\times X^{r_\infty}\srarr X^m\times X^{r_0}\times
X^{r_\infty},$$ and for $i$, an interior marked point, consider
the evaluation map
$$\ev_i:\my\srarr P.$$

Let $\oh(1)$ be the equivariant line bundle over $P$ that is dual
to $\oh(-1)$ equipped with the linearization
$$\begin{array}{rcccl}
& (l,t) & \mapsto & (l,\lambda^{-1}t)&\\
&\downarrow & & \downarrow & \\
& [l:t] & \mapsto & [\lambda l:t]& \\
\end{array}$$
Note that $i_\infty^*(\oh(1)\otimes p^*L)$ is the equivariant
trivial bundle so
$$i_0^*(c_1(\oh(1))+c_1(L))=\hbar+c_1(L)\in \HS{*}(D_0)=H^*(D_0)[\hbar]$$
while
$$i_\infty^*(c_1(\oh(1)+c_1(L))=0 \in
\HS{*}(D_\infty)=H^*(D_\infty)[\hbar].$$

Note also that there is also a natural map
$\pt^*:\com[\hbar]=\HS{*}(\pt)\srarr \HS{*}(\my)$ and that
$\pt^*\hbar$ (which we will denote by $\hbar$) is an equivariant
extension of $0\in H^2(\my)$.

Let $n=\vdim \my$ and $c\in H^{n-4}(X^m\times X^{r_0} \times
X^{r_\infty})$

Note that
$$\deg\left((\hbar \cup \Ev^*c \cup
\ev_i^*(c_1(\oh(1))+c_1(L))) \cap \vir{\my}\right )|_{\hbar=0}=0$$
because the cohomology class is an equivariant extension of $0$.
We will prove Theorem \ref{degenma} by computing this degree by
virtual localization.

\subsection{Fixed Loci}

We identify the $\com^*$ fixed loci in $\my$.  For a map $f$ in a
$\com^*$ fixed locus, we have the composition
$$C\srarr\ _{k,l} Y\srarr .P$$
The irreducible components of $C$ are of three types: (a) those
mapping into $D_0$; (b) those mapping into $D_\infty$; (c) those
mapping into a fiber of $P\srarr X$ totally ramified over two
points.  There are, therefore, three types of fixed loci:
\begin{enumerate}
\item Those whose generic element only has components of type (a)
and (c)

\item Those whose generic elements only has components of type (b)
and (c)

\item Those whose generic element has components of all three
types.
\end{enumerate}

A fixed locus of the first type is parameterized by
$T:\cm(\ca,\Gamma)\srarr F$ where the morphism $T$ attaches
components of type (c) of degrees
$\{\mu^\infty(1),\dots,\mu^\infty(|R_\infty|)\}$ to the boundary
marked points.  We have $T_*(\vir{\ma})=\vir{F}$.  The evaluation
maps fit into the following commutative diagram
$$\xymatrix{
\ma\ar[r] \ar[d] & X^m\times X^{r_0} \times X^{r_\infty} \ar[d]\\
\my\ar[r]        & P^m\times X^{r_0} \times X^{r_\infty} }$$ where
the vertical map $X\srarr P$ is given by including $X$ as
$D_0\subset P$. The virtual normal bundle to this fixed locus has
Euler class
$$e(N)=\hbar-c_1(\Bot).$$

A fixed locus of the second type is parameterized by
$B:\cm(\ca,\Gamma)\srarr F$ where the morphism $B$ attaches
components of degrees $\{\mu^0(1),\dots,\mu^0(|R_0|)\}$. The
evaluation maps fit together as before except that the inclusion
is now $D=D_\infty\hookrightarrow P$. $B_*(\vir{\ma})=\vir{F}$ and
$$e(N)=-\hbar-c_1(\Top).$$

A fixed locus of the third type is parameterized by stacks of the
form $I:\cm(\ca,\Gamma_{A_t})\times_{D^r}
\cm(\ca,\Gamma_{A_b})\srarr F$ corresponding to a graph join
quadruple $(\Gamma_{A_t},\Gamma_{A_b},L,J)$ such that
$\Gamma_{A_t}*_{L,J}\Gamma_{A_b}=\Gamma$. The morphism $I$ inserts
components of type (c) of degrees
$$\{\mu_t^\infty(1),\dots,\mu_t^\infty(|R_{t\infty|})\}=\{\mu_b^0(1),\dots,\mu_b^0(|R_{b0}|)\}$$
between the components coming from each factor of $\ma$. The
evaluation map takes the interior marked points on
$\cm(\ca,\Gamma_{A_t})$ and $\cm(\ca,\Gamma_{A_b})$ to
$X=D_0\subset P$ and $X=D_\infty\subset P$, respectively.  We have
$$I_*(\vir{\cm(\ca,\Gamma_{A_t})\times_{D^r}
\cm(\ca,\Gamma_{A_b})})=|\Aut_{\Gamma_{A_b},\Gamma_{A_t},L}(RA_{b0},RA_{t\infty})|.$$
If $p_t,p_b$ are the projections of
$\cm(\ca,\Gamma_{A_t})\times_{D^r} \cm(\ca,\Gamma_{A_b})$ onto
each factor, then the normal bundle to the fixed locus has Euler
class
$$e(N)=(\hbar-p_t^*c_1(\Bot))(-\hbar-p_b^*c_1(\Top)).$$

\subsection{Localization Computation}

We now compute the contribution from each fixed locus.  The
virtual localization formula \cite{GP} states that given a
top-dimensional class $b\in H^*(\my)$, we have
$$\deg(b\cap\vir{\my})=\sum_{I:F\srarr \my} \frac{1}{\deg(I)}
\deg\left(\frac{I^*b}{e(N_F)}\cap \vir{F}\right).$$

The fixed locus of the first type contributes
\begin{eqnarray*}
&&\deg\left(\frac{T^*(\hbar
\ev_i^*(c_1(\oh(1))+c_1(L))\Ev^*c)}{e(N)} \cap
\vir{\ma})\right)\\
&=&\deg\left(\frac{\hbar(\hbar+\ev_i^*c_1(L)}{\hbar-c_1(\Bot)}
\Ev^*c \cap \vir{\ma})\right)\\
&=&\deg((\ev_i^*c_1(L)+c_1(\Bot))\Ev^*c \cap \vir{\ma}).
\end{eqnarray*}

Fixed loci of the second type do not contribute to the
localization formula because $\ev_i:F\srarr P$ factors as
$$\ev_i:F=\ma\srarr X=D_\infty \hookrightarrow P$$
and $i_\infty^*(c_1(\oh(1))+c_1(L))=0$.

The only fixed loci of the third type that contribute are those in
which the $i$th marked point is mapped to $D_0$. Such a fixed
locus contributes
\begin{eqnarray*}
&&\frac{1}{\Aut}\deg\left(\frac{I^*(\hbar
\ev_i^*(c_1(\oh(1))+c_1(L))\Ev^*c)}{e(N)} \cap
\vir{\ma\times_{D^r}\ma})\right)\\
%%&=&\frac{1}{\Aut}\deg\left(\frac{\hbar(\hbar+\ev_i^*c_1(L))}{(\hbar-p_t^*c_1(\Bot))(-\hbar-p_b^*c_1(\Top))}
%%\Ev^*c \cap \vir{\ma\times_{D^r}\ma})\right)\\
&=&-\frac{1}{\Aut}\deg(Ev^*c\cap \vir{\ma\times_{D^r}\ma})
\end{eqnarray*}
where
$\Aut=|\Aut_{\Gamma_{A_b},\Gamma_{A_t},L}(RA_{b0},RA_{t\infty})|$.

Putting everything together we get
\begin{eqnarray*}
0&=&\deg((\ev_i^*c_1(L)+c_1(\Bot))\Ev^*c \cap
\vir{\ma})\\
&&-\frac{1}{|MA_b|!|MA_t|!(|{RA_{b0}}|)^2}\sum_{\Upsilon\in\Omega_{\Lnb{i}}}\deg(Ev^*c\cap
\vir{\ma\times_{D^r}\ma})
\end{eqnarray*}
where the sum is over all quadruples
$\Upsilon=(\Gamma_b,\Gamma_t,L,J)\in\Omega_{\Lnb{i}}$ as in
Theorem \ref{interpofbundles}.  By choosing $c$ to be a formal sum
of variables as in the definition of the correlators, we get
Theorem \ref{degenma}.

%\end{document}